\documentclass[11pt]{article}

\usepackage{graphicx}
\usepackage[margin=1in]{geometry}
\usepackage{amsmath,amssymb,amsthm,mathtools}
\usepackage{threeparttable}
\usepackage{xcolor}
\usepackage{enumitem}
\usepackage[hidelinks]{hyperref}
\usepackage{bbm}
\usepackage{booktabs}
\usepackage{tabularx}

\newtheorem{theorem}{Theorem}
\newtheorem{lemma}{Lemma}

\newtheorem{corollary}{Corollary}
\newtheorem{conjecture}{Conjecture}
\theoremstyle{definition}
\newtheorem{definition}{Definition}
\theoremstyle{remark}
\newtheorem{remark}{Remark}
\newtheorem{example}{Example}

\DeclareMathOperator{\diag}{diag}
\DeclareMathOperator{\Var}{Var}
\newcommand{\one}{\mathbbm{1}}          % vector of ones (was e)
\newcommand{\R}{\mathbb R}
\newcommand{\E}{\mathbb E}
\newcommand{\Prob}{\mathbb P}
\newcommand{\dt}{\Delta t}
\newcommand{\sspc}{\mathcal C}          % SSP coefficient
 
\newcommand{\dtFE}{\Delta t_\textup{FE}}

\newcommand{\leanrepo}{https://github.com/RKProofs/Lean-Files}
\newcommand{\leandir}{SSPLean}
\newcommand{\leandocdir}{doc}
\newcommand{\leanbranch}{arXiv}
\newcommand{\leanproof}[3]{\href{\leanrepo/blob/\leanbranch/\leandir/#2.lean\#L#3}{\texttt{#1}}}
\newcommand{\leandocfile}[1]{\href{\leanrepo/blob/\leanbranch/\leandocdir/#1}{\texttt{#1}}}
\newcommand{\leanreporef}{\href{\leanrepo/blob/\leanbranch}{\texttt{\leanrepo}}}

\newcommand{\resref}[2]{\hyperref[#2]{#1~\ref*{#2}}}

\title{A stability order barrier and SSP bound for Runge--Kutta methods
with non-negative coefficients: AI-assisted proofs}
\author{Sigal Gottlieb\thanks{Department of Mathematics, UMass Dartmouth} , 
Zachary Grant$^*$, Lior Horesh\thanks{IBM Research,
Yorktown Heights, New York} \\
Lajos L\'oczi\thanks{CEMSE division,KAUST} , 
Avraham Shinnar$^\dagger$, Barry Trager$^\dagger$}
\date{\today}

\begin{document}

\maketitle

\begin{abstract}
We prove three results about Runge--Kutta methods with non-negative coefficients.
First, every $s$-stage implicit Runge--Kutta method of order at least two has SSP
coefficient at most $2s$; the bound is sharp and is attained by a known diagonally
implicit method, so full implicitness gives no advantage over DIRK at second order.
This extends the second-order DIRK result of Kocsis \& N\'emeth to full coefficient
matrices and settles a conjecture of Ketcheson, Macdonald \& Gottlieb.
These proofs  were found with the assistance of AI tools, which uncovered an
interesting relationship between the Butcher matrices for SSP methods and 
phase-type distributions, and led to the proof of the second result: 
a bound on the SSP coefficient of $s$-stage implicit Runge--Kutta methods of 
order at least three. 
The third result in this work is a proof that a third-order Runge--Kutta method
with entrywise non-negative $A$ and $b$ cannot be  algebraically stable.
This algebraically intensive proof was also AI-assisted. 
This paper  describes the iterative process of AI-assisted results, 
and presents an approach  that centers  human understanding and creativity in the  
AI-assisted mathematical endeavor.

\end{abstract}

%\david{Any reason why we are not using bibtex?  In my opinion it will save a lot of time and it would be easier to switch now rather than later.}
%\sigal{Certainly we can use bibtex. Would you be willing to put it in?}

\begin{quote} {\em
    This paper is dedicated to the memory of Marc Spijker (1939--2026) with appreciation
    for his work on strong stability preserving time stepping methods.}
\end{quote}

\section*{Acknowledgments}
The authors gratefully acknowledge David Ketcheson, who was part of the collaboration from its inception, provided substantive intellectual input, and participated in writing the manuscript. His insightful suggestions shaped both the content and development of this project, and his contributions were essential 
to the paper.

SG and ZJG's research was supported in part by  AFOSR Grants No. FA9550-23-1-0037
and No. FA9550-26-1-B182, NSF Grant No.  DMS-2309609, 
and DOE  Grant No. DE-SC0023164 Subaward RC114586.

\section{Introduction}
The focus of this  paper is on the use of AI tools to prove 
conjectures on Runge--Kutta methods that have non-negative coefficients.
Such methods may have strong stability preserving (SSP) properties 
and so have proven useful in the time evolution of hyperbolic partial
differential equations (PDEs).  The first set of proofs make use of a 
dictionary between SSP theory   and absorbing Markov chains, that was revealed by
AI tools. In this dictionary (described in Section \ref{sec:dictionary}),
the SSP conditions assert that a related matrix is substochastic, 
the order conditions become moment conditions on a phase-type absorption time, 
and the extremal methods correspond to Erlang distributions.  
The key inequality introduced by this dictionary is the Aldous--Shepp 
least-variable-phase-type theorem.  
These connections to  substochastic matrices allow us to prove a
\emph{hierarchy} for acyclic phase-type moment inequalities whose first member
is the Aldous--Shepp theorem \cite{AldousShepp1987}, 
and show that its second member implies the sharp bound
$r \le s-1+\sqrt{s^2-1}$ on diagonally implicit Runge--Kutta  methods.
We also show that the second member is true for the cyclic case, 
which proves the sharp bound previously conjectured for third-order 
implicit Runge--Kutta methods. The use of the connections revealed 
by AI tools and described in Section \ref{sec:dictionary} allowed us to extend our results to third order, which we then found matched our past numerical results perfectly. 
The next result, in Section \ref{sec:barrier},  explores the failure of
numerical searches to find third order algebraically stable 
implicit Runge--Kutta methods with  non-negative coefficients,
and prove that the bushy tree conditions for order three
force a conflict between SSP and algebraic stability.

The elementary parts of all proofs, including the order 
barrier and the reductions from moment inequalities to SSP bounds, 
have been machine-verified in Lean~4.
%the formalization process itself exposed and repaired several gaps in earlier drafts.
We used mainly ChatGPT, Claude, and Gemini, augmented by Lean, 
to explore the problems, generate candidate arguments, and verify the results
by human and Lean formalization.  
In our exposition, we have sought to  ``increase the emphasis we place on proof digestion'' as advised by \cite{Tao2026AI}, and we
hope that our presentation of the results as well as the process, 
``enables people to understand and think more clearly and effectively about
math'' \cite{Thurston}. 
This paper presents more than the proofs: it  presents  our approach to utilizing 
AI tools in a rigorous and thoughtful way that enhances human creativity and 
develops mathematical connections  and insight  as well as results.

\subsection{On AI as part of the evolution of computational tools}
It is interesting to note the evolution of computation in the 
field of SSP methods. In 1992, SG first used an optimization code 
to try to find a SSP decomposition of  the four stage fourth order 
explicit Runge--Kutta method. This of course failed, but pointed
the way to a proof that there are no SSP four stage fourth order explicit
Runge--Kutta methods. Later, Ferracina and Spijker \cite{FerracinaSpijker2004,FerracinaSpijker2008} 
reformulated the optimization problem, and DK refined this approach and 
generated an optimization suite for SSP methods of many types \cite{RKopt}.
In fact, the central conjecture of this paper, the ``$2s$ conjecture''
was first proposed by Ferracina and Spijker for 
singly diagonally implicit Runge--Kutta (SDIRK) methods 
as a result of their  numerical searches, and later expanded to encompass 
DIRK and even fully implicit methods
by \cite{KMG09} as a result of their numerical searches. 
This paper represents two advances in the use of computation: 
the inclusion of AI tools to find connections and
proofs of facts about Runge--Kutta methods with non-negative coefficients, 
and the use of Lean to formalize these.

The proofs written by AI were rarely readable or intuitive;
we digested the proofs and ensured that each result was checked, 
corrected, formalized, written,  re-written,  understood, and explained
by a human.  To the best of our ability, we tried to present results 
that were interesting and that revealed new connections and new mathematical results. 
We studied the significance of the observations 
and connections uncovered by AI and used the insights we 
gained from this process to deepen our understanding of the field.
In all the above, the ongoing active involvement of several 
mathematicians  enabled rigorous  and comprehensible mathematics, 
which was not produced by any of the AI tools.
Our approach centers  human understanding and creativity in the 
mathematical endeavor, with the  AI tools serving a valuable assistive role in 
revealing connections and in handling complicated algebraic arguments.
We believe that this approach is a valuable one.

\section{Runge--Kutta methods and strong stability preservation}
\label{sec:background}

An $s$-stage Runge--Kutta method applied to $u'(t) = f(u(t))$ can be written as
\begin{align} \label{eq:RK}
    u^{(i)} & =  u_n + \dt \sum_{j=1}^{s} a_{ij} f(u^{(j)}), \qquad i = 1,\dots,s,
    \notag\\
    u_{n+1} & =  u_n + \dt \sum_{i=1}^s b_{i} f(u^{(i)}) .
\end{align}
The coefficients form a matrix $A = (a_{ij})$ and a vector $b = (b_i)$; we write
$\one$ for the vector of ones, $c = A\one$ for the abscissas, and all vector and
matrix inequalities are understood componentwise.  The order conditions through
second order are
\begin{equation} \label{eq:OC}
    b^\top \one = 1, \qquad
    b^\top c = \tfrac12 . \qquad
%    b^T A c = \tfrac16, \qquad
%    b^T (c \circ c)  = \tfrac13 ,
\end{equation}
%where $\circ$ is the componentwise product.  
We are interested in methods whose
coefficients are non-negative: $A \ge 0$ and $b \ge 0$.

A Runge--Kutta method with Butcher coefficients $(A,b)$ may contain redundant stages that can be removed or merged without changing its numerical one-step map. Two standard notions of reducibility describe these possibilities \cite{HairerWanner}. HS-reducibility concerns stages that are equivalent and can be merged. Identical rows of $A$ provide a simple example, but stages may also be equivalent when their coefficient rows become identical after coefficients associated with equivalent stages are grouped and summed. Linear dependence among rows of $A$, by itself, does not imply that a stage can be eliminated, since its function evaluation may still be needed. DJ-reducibility concerns stages that do not influence the final output, either directly through their weights or indirectly through other stages that contribute to the output.
We use ``irreducible'' to mean that neither HS-reduction nor DJ-reduction is possible. This distinction is particularly relevant when studying bounds that depend on the number of stages: redundant stages increase the stated stage count without changing the numerical method. In the remainder of the paper, whenever we discuss a Runge--Kutta method, 
we assume that it is irreducible.

\subsection{Necessary and sufficient conditions for SSP methods}
Strong stability preserving Runge–Kutta methods were developed  
\cite{ShuOsher1988,ShuOsher1989}
to preserve the nonlinear stability properties of forward Euler (in any convex functional)
under higher order time stepping. 
This approach was deeply studied for Runge–Kutta, multistep, general linear methods, 
and multi-derivative methods.
%More recently, as well as other extensions including integrating factor methods and linear splitting methods. 
These methods have proven useful in the simulation of hyperbolic PDEs in
conjunction with  many spatial discretizations
and are widely used in many application areas \cite{GKSbook}.  

\begin{definition}
The SSP coefficient $\sspc(A,b) $
(or \emph{radius of absolute monotonicity}) of a Runge--Kutta method is the largest value $r\ge 0$ such that
\begin{subequations}\label{eq:KFSconditions}
\begin{eqnarray}
    (I + rA) \; \; \; \mbox{is invertible} \\
    (I + rA)^{-1} \one \geq  0 \\
    A (I + rA)^{-1}  \geq 0 \\
    b^\top  (I + rA)^{-1} \geq 0 \\
    1- r b^\top  (I + rA)^{-1} \one \geq 0,
\end{eqnarray}
\end{subequations}
where the inequalities are understood componentwise.
\end{definition}

The SSP coefficient is relevant for initial value problems (IVPs)
\begin{align} \label{ivp}
    u'(t) & = f(u) \quad 
    u(t_0)  = u_0 
\end{align}
for which $f$ satisfies a monotonicity condition
(or \emph{circle condition})
in some norm or semi-norm $\|\cdot\|$:
\begin{align} \label{FEcond}
    \|u + \dt f(u)\| \le \|u\|  \quad \quad \text{ for all } \dt \in [0, \dtFE].
\end{align}
Here the notation $\dtFE$ is motivated by the fact that \eqref{FEcond} is monotonicity of a forward Euler step.
We define ${\mathcal F}_{\dtFE}$ as the set of all IVPs satisfying \eqref{FEcond} for a given maximal step size $\dtFE$.

\begin{theorem}\cite{Kraaijevanger1991, FerracinaSpijker2004}\label{sspthm}
Let an irreducible Runge--Kutta method defined by $(A,b)$,
and a function $f \in {\mathcal F}_{\dtFE}$ be given.
Then the numerical solution obtained by applying the method $(A,b)$ to the initial-value problem \eqref{ivp}
satisfies
 \begin{subequations} \label{monotonicity}
 \begin{align}
    \|u^{(i)}\| & \le \|u_n\|  \label{Monotonicity:Internal}\\
     \|u_{n+1}\| & \le \|u_n\| \label{Monotonicity:Step}
 \end{align}
 \end{subequations}
 if the step size satisfies
 \begin{align} \label{step-size-restriction}
     0 \le \dt \le  \sspc(A,b)  \; \dtFE.
 \end{align}
Moreover, this step-size coefficient is sharp 
in the sense that for every finite $\theta > \mathcal C(A,b)$, 
there exist a space, a seminorm, 
a function satisfying Condition \ref{FEcond}, 
and initial data for which the Runge–Kutta equations \eqref{eq:RK} admit a solution 
at $\Delta t=\theta\Delta t_{\mathrm{FE}}$ that violates at least one of the monotonicity conditions \eqref{monotonicity}.
\end{theorem}

\subsection{The $\sspc(A,b)  \leq 2s$ conjecture}
% \david{Here we introduce $\sspLin$ as a dummy variable and $\rho$ as the threshold factor.  So I started cleaning up the text to be consistent with that.  But then it seems like through most of the paper we use $\sspLin$ as the threshold factor.  Please choose a notation and use it consistently.}

Theorem \ref{sspthm} says that the Euler condition \eqref{FEcond} and step size restriction \eqref{step-size-restriction} imply monotonicity \eqref{monotonicity} for all $f$.  To describe the origin of the conjecture, we consider two simplifications.
The first simplification occurs when we restrict our attention to linear, 
autonomous functions $f$ and  consider the 
\emph{threshold factor}:
\begin{align}
    \rho(A,b) = \sup\{r \ge 0 : 0 \le \dt \le r \dtFE \implies \text{ \eqref{Monotonicity:Step} for all linear, autonomous IVPs\}}
\end{align}
%largest value of $\sspLin$  such that \eqref{FEcond} and 
%$$  0 \le \dt \le \sspLin \dtFE $$
%imply step-wise monotonicity \eqref{Monotonicity:Step}
%over this smaller class of problems.
%In this case, the upper bound on $\sspLin$ is given by 
%what is known as the \emph{threshold factor}, denoted
%herein by $\rho$.
By definition, $\rho$ is also an upper bound on $\sspc(A,b)$.

For explicit methods, it is straightforward to show that $\rho$ (hence $\sspc$)  is bounded by $s$, the number of stages of the method.  
Since the explicit Euler method has (by definition) $\rho=1$, 
and since the cost of the method scales linearly with $s$, explicit methods can at most hope to achieve the step size efficiency of Euler's method in this respect.  
Nevertheless, high order explicit SSP methods are extremely useful since 
the first-order explicit Euler method is too inaccurate for practical integration.
On the other hand, the implicit Euler method has $\rho=\infty$, and indeed
$\sspc = \infty$. For general (implicit) methods, it was shown early 
on that for methods of order 2 (and higher), $\rho < \infty$; see \cite{bolley1978,spijker1983}.  
In 1986 it was conjectured that $\rho\le 2s$ for implicit methods \cite{vandegriend1986}; this conjecture was later disproved \cite{2014_rational_absmon}, but available counterexamples exceed $2s$ by only a small amount. 

Meanwhile, numerical optimization of  $\sspc$ for high-order 
implicit Runge--Kutta methods yielded no methods with $\sspc>2s$, 
leading to the new conjecture that $\sspc \le 2s$ for all methods of order $p\ge2$.

\begin{conjecture}[The $2s$ conjecture] \label{conj:SSPbound}
    Given a scalar $r > 0$ and 
    matrices $A \in \mathbb{R}^{s\times s}$ and  $b \in \mathbb{R}^{s}$ 
that satisfy the SSP conditions \eqref{eq:KFSconditions}:
\[   (I + rA) \; \; \; \mbox{is invertible} ,\]
 \[   (I + rA)^{-1} \one \geq  0 ,\quad
    A (I + rA)^{-1}  \geq 0 ,\quad
    b^\top  (I + rA)^{-1} \geq 0, \quad
    1- r b^\top  (I + rA)^{-1} \one \geq 0, \]
(where the inequalities are understood componentwise),
and the order conditions 
\[b^\top  \one = 1 , \qquad   b^\top  A \one  =   \frac{1}{2}  \]
 then $r \leq 2s$.
\end{conjecture} 

Curiously, existing numerically optimized implicit SSP methods all have 
the (desirable)  property of being diagonally implicit; i.e. 
the coefficient matrix $A$ is lower triangular.  
Furthermore, in numerical optimization tests we found that this
bound holds not only for all known Runge--Kutta methods, 
but also for all known general linear methods.  
The $2s$ bound represents a strong fundamental limitation on 
numerical integration methods.  
If one requires strict monotonicity with respect to an {\em inner product norm}, 
there exist methods of arbitrarily high order that allow arbitrarily large step sizes.  
But for more general strict boundedness properties like \eqref{monotonicity} in non-inner-product norms, the efficiency of higher-order implicit methods is extremely restricted.

\subsection{Proof of the SSP bound $r \le 2s$ for implicit methods}
\label{sec:SSPbound}

The first main result of this work
%given in Section \ref{sec:SSPbound}, 
is a proof of the $2s$ conjecture \ref{conj:SSPbound}.
Substantial progress on the conjecture was made by
Kocsis \& N\'emeth, who proved that it holds for the class of diagonally implicit (DIRK) methods \cite{KocsisNemeth}.
In other words, whereas Kocsis \& N\'emeth require the matrix $A$ to be lower triangular,
we are able to remove this restriction and extend the result to methods with a full matrix $A$. As much as  possible, we follow their notation and process. 

As in \cite{KocsisNemeth}, we translate the  conditions  
\eqref{eq:KFSconditions} for any $r>0$ 
to conditions on the  invertible matrix $N\in \mathbb{R}^{s\times s}$ and
vector $w \in \mathbb{R}^{s}$ 
defined by:
\begin{equation}\label{eq:NQw}
  N = I + rA, \qquad  w^\top  = r\, b^\top  N^{-1}.
\end{equation}
Noting that $I - N^{-1} = rA(I+rA)^{-1}$, the absolute monotonicity conditions \eqref{eq:KFSconditions} become:
\begin{subequations} \label{NemethSSPconditions}
\begin{eqnarray} 
    N \; \; \; \mbox{is invertible}  \label{NSSP1}\\
    N^{-1} \one \geq  0 \label{NSSP2} \\
%        (N - I) N^{-1}  \geq 0 
    I -  N^{-1}  \geq 0  \label{NSSP3} \\
    w^\top  \geq 0  \label{NSSP4}\\
    1- w^\top  \one \geq 0. \label{NSSP5}
\end{eqnarray}
\end{subequations}
Again, the inequalities are understood componentwise.
The first order condition  $b^\top  \one = 1 $  can be rephrased as
\[w^\top  N \one = r b^\top  N^{-1} N \one = r b^\top  \one = r,\]
while the second order condition
$ b^\top  c = \tfrac12$ is equivalent to
\[ w^\top  N^2 \one =  r b^\top  N^{-1} N^2 \one = r b^\top  N \one = r b^\top  \left( I + rA\right) \one 
= r b^\top  \one + r b^\top  r A \one  = r + r^2 b^\top  c  = r + \frac{1}{2} r^2.
\]
So the order conditions up to second order are equivalent to
\begin{subequations} \label{NemethOrderConditions}
\begin{eqnarray} 
w^\top  N \one & = &  r \\
w^\top  N^2 \one & = &  r + \frac{1}{2} r^2.
\end{eqnarray}
\end{subequations}

We now state our first main result: this is a restatement of Conjecture \ref{conj:SSPbound} above in the transformed variables. It is also a generalization of Theorem 2.6 in 
Kocsis \& N\'emeth \cite{KocsisNemeth} to the fully implicit case.

\begin{theorem} \label{thm:SSPbound}
Given matrices $N \in \mathbb{R}^{s\times s}$ and  $w \in \mathbb{R}^{s}$ 
defined by \[ N= I + rA, \; \; \;  w^\top  = r b^\top  N^{-1}\] 
that satisfy the SSP conditions \eqref{NemethSSPconditions}
\[    N \; \; \; \mbox{is invertible} , \qquad
    N^{-1} \one \geq  0 ,  \qquad
    Q = I -  N^{-1}  \geq 0  , \]
\[    w^\top  \geq 0 , \qquad   1- w^\top  \one \geq 0, \]
(where the inequalities are understood componentwise),
and the order conditions \eqref{NemethOrderConditions}
\[w^\top  N \one = r , \qquad   w^\top  N^2 \one  =  r + \frac{1}{2} r^2, \]
 then $r \leq 2s$.
\end{theorem}

In this  section we assemble the elements needed to prove  Theorem \ref{thm:SSPbound}. 
A key component of this proof is the central bound, the substochastic sequence inequality (Lemma \ref{lem:AS}), that uses  
uses an important  connection between SSP conditions and absorbing Markov chains.
This connection was not previously known in the SSP/monotonicity field and 
was revealed by ChatGPT.
We present each of the results leading up to the proof separately,
to show which results use this newly revealed  connection to absorbing 
Markov chains, and which do not. 
The aim in this section is to use the connections as minimally as we can, 
and stick closely to the notation and concepts more usual in SSP theory.
In Section \ref{sec:dictionary}  we describe the dictionary between 
SSP conditions and absorbing Markov  chains,
and rely more heavily on these connections to prove this result 
more elegantly and with some more generality. 
%\revnote{Is it more generality? Check this --SG}

The following  lemma  is closely related to \cite[Lemma 2.9]{KocsisNemeth}.  
It is slightly simpler, because it
avoids using any expressions in which the denominator could vanish.
%We state it in a division-free form: the ratio form
%$(w^\top  N \one)^2/(w^\top  N^2 \one) \le \max_i \big((N \one)_i^2/(N^2 \one)_i\big)$ 
%requires additional care when denominators vanish (e.g.\ $w = 0$), 
%whereas the form below is more robust, and still well-suited to the 
%proof of Theorem \ref{thm:SSPbound}.

\begin{lemma}[weighted Cauchy--Schwarz bound]\label{lem:wcs}
Let $u, v, w \in \R^s$ with $w \ge 0$, $w^\top \one \le 1$, and suppose
$u_i^2 \le K v_i$ for all $i$ and some $K \ge 0$.  Then
\[ (w^\top  u)^2 \;\le\; K\, w^\top  v . \]
\end{lemma}

\begin{proof}
By Cauchy--Schwarz with the splitting $w_iu_i = \sqrt{w_i}\;(\sqrt{w_i}\,u_i)$,
\[
(w^\top u)^2 = \Bigl( \sum_i \sqrt{w_i} (\sqrt{w_i}\,u_i) \Bigr)^2 
\le \Bigl(\sum_i w_i\Bigr)\Bigl(\sum_i w_i u_i^2\Bigr)
\le \sum_i w_i u_i^2
\le K \sum_i w_i v_i  =  K\, w^\top  v ,
\]
using $w^\top \one \le 1$ together with $\sum_i w_iu_i^2 \ge 0$ for the middle step, and
$w_i \ge 0$, $u_i^2 \le Kv_i$ for the last.
\end{proof}

We now state some lemmas that elucidate the behavior of the  relevant matrices.
\begin{lemma}[substochastic structure]\label{lem:substochastic}
Suppose $N = I + rA$ satisfies \eqref{NSSP1}--\eqref{NSSP3} for some $r>0$, and
let $Q = I - N^{-1}$.  Then 
$Q$ satisfies the componentwise inequalities
\[ Q \one \leq \one \; \; \; \mbox{and} \; \; \;  0 \leq Q_{ij} \leq 1 . \]
Furthermore, the spectral radius of $Q$  is bounded:
$ \rho(Q) < 1 ,$
so that 
\[    N = (I-Q)^{-1} = \sum_{k=0}^{\infty} Q^k \;\ge\; 0 . \]
\end{lemma}

\begin{proof} 
Using the property $ N^{-1} \one \ge 0 $ and the definition of $N= (I - Q)^{-1}$ we have 
\[ (I - Q) \one  \geq 0 \; \; \Rightarrow \; \; Q \one \leq \one . \] 
The condition  \eqref{NSSP3} is exactly $Q \geq 0$. 
Since each $Q_{ij}$ is non-negative we have
\[ Q_{ij}  \leq \sum_{j=1}^s Q_{ij} = (Q \one)_i \leq \one_i = 1 \; \; \;
\Rightarrow  \; \; \;  0 \leq Q_{ij} \leq 1.\]

These conditions on the elements of $Q$ give us a bound on its maximum norm:
$ \|Q\|_\infty = \max_{1\le i\le s}\sum_{j=1}^s |Q_{ij}| =
\max_{1\le i\le s}\sum_{j=1}^s Q_{ij} \le 1. $
The spectral radius is bounded by the maximum norm  so that 
\[ \rho(Q)\le \|Q\|_\infty\le 1. \]
The Perron–Frobenius theorem states that for a matrix $Q$ with elements
$Q_{ij} \geq 0$ 
the spectral radius $\rho(Q)$ is a real non-negative eigenvalue of $Q$, and  
in particular, there exists a nonzero vector $v\ge 0$ such that
$Q v=\rho(Q) v$.  Assume that $\lambda =1 $ is an eigenvalue of $Q$
with an associated nonzero eigenvector $v_\lambda$ so that
then $Q v_\lambda = v_\lambda $. But this implies that 
$(I- Q) v_\lambda = 0 $ where $v_\lambda$ is nonzero and
$N^{-1}  = I-Q$ is invertible. This is a contradiction 
so we must have $\rho(Q) < 1.$
\end{proof}

\begin{lemma}[Neumann series for the square]\label{lem:neumann2}
If $Q\in\mathbb{R}^{s\times s}$ satisfies
\[ Q_{ij} \ge 0, \qquad Q\one \le \one\; \; \; \mbox{and } \; (I - Q) \;  \mbox{is invertible}, \]
then $\rho(Q)<1$ and
\[ (I-Q)^{-2}= \sum_{k=0}^\infty (k+1)\, Q^k. \]
\end{lemma}

\begin{proof}
The proof of Lemma \ref{lem:substochastic} gave us $\rho(Q)<1$.  
For $|t| < 1/\rho(Q)$ define the absolutely convergent Neumann series
$F(t)=(I-t Q )^{-1} = \sum_{k\ge0} t^k Q^k,$
and differentiate it term-by-term (which is valid inside the radius of
convergence) to obtain
\[ F'(t)  = \sum_{k=1}^\infty k t^{k-1}Q^k,\]
and also differentiating $F(t)$ directly
\[ F'(t)  = Q(I-tQ)^{-2}, \]
so that 
\[  Q(I-tQ)^{-2} = \sum_{k=1}^\infty k t^{k-1}Q^k. \]
Now observe that
%\begin{eqnarray*}
\[ (I- t Q)^{-2}   = t Q (I - t Q)^{-2}  + (I- t Q)^{-1}  
=  \sum_{k=0}^\infty k\,t^{k} Q^k + \sum_{k=0}^\infty t^k Q^k    
=  \sum_{k=0}^\infty (k+1) t^k Q^k.
\]
Finally, since $\rho(Q) < 1$ we may set  $t=1$ to obtain the result.
We note that the last series expansion is simply a consequence of 
$\rho(Q)<1$.

A simpler proof, but one that requires more sophisticated counting, is given by 
    \[ (I-Q)^{-2}= \left(\sum_{k=0}^\infty \, Q^k\right)^2
    =  \sum_{j=0}^\infty \sum_{k=0}^\infty \, Q^{j+k}
    = \sum_{k=0}^\infty (k+1)\, Q^k .\]
\end{proof}

The following lemma allows us to bound the needed elements.
The proof of this lemma contains the key element provided by  
ChatGPT: the translation of our problem 
to a different mathematical language that enables us to use mathematical tools that are not typically 
seen in this context. 
\begin{lemma}[substochastic sequence inequality]\label{lem:AS}
Let $Q$ be an $s \times s$ matrix with $Q \ge 0$, $Q\one\le \one$ and $\rho(Q) < 1$.
Fix an index $1\le i\le s$ and set $p_k=\bigl(Q^{k}\one\bigr)_i$ for $k \ge 0$.  Then
\[ \left(\sum_{k=0}^{\infty}p_k\right)^{2}
\;\le\; \frac{2s}{s+1}\sum_{k=0}^{\infty}(k+1)\,p_k . \]
\end{lemma}

\begin{proof}
Let $G=Q-I$.  The off-diagonal entries of $G$ are non-negative and $G\one\le0$, so
$G$ is the transient generator of a continuous-time Markov chain on $\{1,\dots,s\}$
with one additional absorbing state \cite[Ch.~3]{KemenySnell1976}; the absorption
rate from state $j$ is $1-\sum_{\ell} q_{j\ell}\ge0$.  Since $\rho(Q)<1$, every
eigenvalue of $G$ has negative real part, so the absorption time is almost surely
finite with finite moments of every order.

Start the chain in state $i$ and let $T$ be its absorption time; 
writing $\delta_i$
for the $i^\text{th}$ coordinate vector, the survival probability is
$\Prob_i(T>t)=\delta_i^\top \exp(Gt)\, \one$.  Then

%\lajos{the notation in Sections \ref{sec:dictionary}--\ref{sec:altproof} seems to be different, e.g. no $[\cdot]$ is used after $\E$, I'd suggest using $[\cdot]$}
% \sigal{I now made this change everywhere (I think)}

%\lajos{am I right that $(\cdot)_i$ is the $i^\text{th}$ component of a vector, except for $\delta_i$?}
%\sigal{here everything is with the language of the markov chain, so we
%need not have vectors, just expectations and states}

\[
    \E_i[T]
    =\int_0^\infty \Prob_i(T>t)\,dt
    =\delta_i^\top(I-Q)^{-1}\one
    =\sum_{k=0}^{\infty}(Q^k\one)_i
    =\sum_{k=0}^{\infty}p_k ,
\]
using the Neumann series, and similarly, using Lemma \ref{lem:neumann2},
\[
    \E_i[T^2]
    =2\int_0^\infty t\,\Prob_i(T>t)\,dt
    =2\,\delta_i^\top(I-Q)^{-2}\one
    =2\sum_{k=0}^{\infty}(k+1)p_k .
\]
Note $\E_i[T] = (N\one)_i \ge (\one)_i = 1 > 0$.  The Aldous--Shepp
least-variability theorem \cite{AldousShepp1987} (restated as Theorem \ref{thm:AldousShepp} below) states that the absorption time of
a continuous-time Markov chain with at most $s$ transient states satisfies
\[ \Var_i(T) \ge (\E_i[T])^2/s;\]
hence
\[ \E_i[T^2] \ge \frac{s+1}{s}\,(\E_i[T])^2, \] which is the claim.
\end{proof}
\begin{remark}
The proof of Lemma \ref{lem:AS} depends strongly on 
absorbing Markov chains and related probability theory. The translation
from the SSP language to that of absorbing Markov chains is described below
in Section \ref{sec:dictionary}. This lemma, as most others in this paper,
is formally verified by Lean \cite{lean4} (see Section \ref{sec:lean}).
\end{remark}

\begin{lemma}[componentwise bound]\label{lem:componentbound}
If $Q \ge 0$, $Q\one \le \one$ and $I - Q$ is invertible, then
$N = (I-Q)^{-1}$ satisfies, for every $i$,
\begin{equation}\label{eq:componentbound}
  \bigl((N \one)_i\bigr)^2 \;\le\; \frac{2s}{s+1}\, (N^2 \one)_i .
\end{equation}
\end{lemma}

\begin{proof}
By Lemmas \ref{lem:substochastic}--\ref{lem:neumann2},
$(N\one)_i = \sum_k p_k$ and $(N^2\one)_i = \sum_k (k+1)p_k$ with
$p_k = (Q^k\one)_i \ge 0$; apply Lemma \ref{lem:AS}.
\end{proof}

\begin{remark}\label{rem:dictionaryAldousSchepp}
The inequality of Lemma \ref{lem:componentbound} interpreted 
with the dictionary in Section \ref{sec:dictionary} is a consequence 
of the Aldous--Shepp theorem (Theorem \ref{thm:AldousShepp}). 
Furthermore, we have a Lean proof 
%\sigal{Barry and Avi -- put in citation or link here} 
from Claude that Lemma 5 is in fact equivalent to Aldous--Shepp (\leanproof{aldous\_shepp\_iff\_componentBound\_at}{AldousSheppConverse}{97}).
% The formalization in Lean incidentally yielded an
% elementary proof of the Aldous--Shepp  theorem from 1987 
% without any reference to martingale theory 
% \barry{yes, after a day of work }
We note that more recently we were able with Gemini's help
to prove  Lemma \ref{lem:componentbound} without the new language; 
we present this proof in Appendix \ref{sec:AppendixA}. 
% \lajos{as I see now that proof has been simplified by using ``... so that the assumption
% that $c_1\le c_2\le\cdots\le c_s$ ...'', which inequality is first mentioned only in Remark \ref{rem:cjnondecreasing}, but that Remark is not in this Section, so the reader won't understand the assumption if they jump to Appendix \ref{sec:AppendixA} from here}
%\sigal{I added a comment in the appendix that explains the assumption and points the reader to the remark.}
\end{remark}

\noindent Finally, we are ready to prove Theorem \ref{thm:SSPbound}. 

%\begin{theorem*}\label{thm:main}

%\end{theorem*}

\begin{proof} [Proof of Theorem \ref{thm:SSPbound}]
Let $Q = I - N^{-1}$; by Lemma \ref{lem:substochastic}, $N = \sum_k Q^k \ge 0$, so
$N\one \ge 0$ and $N^2\one \ge 0$, and $r = w^\top N\one \ge 0$.  If $r = 0$ we are
done, so assume $r > 0$.  By Lemma \ref{lem:componentbound}, the vectors
$u = N\one$, $v = N^2\one$ satisfy $u_i^2 \le K v_i$ with $K = \frac{2s}{s+1}$, so
Lemma \ref{lem:wcs} gives
\[
  r^2 = (w^\top N\one)^2 \;\le\; \frac{2s}{s+1}\; w^\top N^2\one
      = \frac{2s}{s+1}\Bigl(r + \frac{r^2}{2}\Bigr).
\]
Multiplying by $(s+1)/r > 0$ and simplifying: $(s+1)r \le 2s + sr$, i.e.\ $r\le 2s$.
\end{proof}

\begin{corollary}[unconditional bound]\label{cor:unconditional}
Every irreducible $s$-stage Runge--Kutta method $(A,b)$ of order $p \ge 2$ has
\[ \sspc(A,b) \;\le\; 2s . \]
\end{corollary}
%\lajos{$\sspc$ has not been defined yet}
% \sigal{I have now included it into the original definition}

\begin{proof}
Let $r \in (0, \sspc(A,b))$.
By definition, conditions \eqref{eq:KFSconditions} 
hold at $r$.  With $N, Q, w$ as in \eqref{eq:NQw}, conditions 
\eqref{eq:KFSconditions} become conditions \eqref{NemethSSPconditions},
and the order conditions become \eqref{NemethOrderConditions}.
Thus, the hypotheses of Theorem \ref{thm:SSPbound} hold.
 Hence $r \le 2s$ for every
$r < \sspc(A,b)$, so $\sspc(A,b) \le 2s$.
\end{proof}

% \david{The following is NOT a corollary; it is a pertinent example.  It has appeared already in the literature, more than once.  If you agree, I would like to relabel it as an "Example" and add citations.  I don't think it's necessary to give the Butcher coefficients, since they are readily available in exactly the same form from other sources.  We should also mention that this method is 2nd order and that it is just a well-known 1-stage method repeated $s$ times, with a scaled step size.}

\noindent The following example shows that this bound is sharp:

\begin{example}[sharpness]\label{cor:sharp}
The bound is attained by the $s$-stage second-order method given in Shu-Osher form
as \cite{KMG09}
\begin{align*}
u^{(1)}
&= u^n + \frac{\Delta t}{2s} f(u^{(1)}), \\[3pt]
u^{(i)} &= u^{(i-1)} + \frac{\Delta t}{2s}
 \left(f(u^{(i-1)})+f(u^{(i)})\right),
 \qquad i=2,\ldots,s, \\[3pt]
u^{n+1} &= u^{(s)}+\frac{\Delta t}{2s} f(u^{(s)}).
\end{align*}
The SSP coefficient is clearly
\[
\mathcal{C}(A,b) = 2s.
\]
% The Butcher coefficients of this method are given by
% \[
% A=
% \begin{pmatrix}
% \dfrac{1}{2s} & 0 & 0 & \cdots & 0 \\[6pt]
% \dfrac{1}{s} & \dfrac{1}{2s} & 0 & \cdots & 0 \\[6pt]
% \dfrac{1}{s} & \dfrac{1}{s} & \dfrac{1}{2s} & \ddots & \vdots \\[6pt]
% \vdots & \vdots & \ddots & \ddots & 0 \\[6pt]
% \dfrac{1}{s} & \dfrac{1}{s} & \cdots & \dfrac{1}{s} &
% \dfrac{1}{2s}
% \end{pmatrix},
% \qquad
% \mathbf{b}
% =
% \begin{pmatrix}
% \dfrac{1}{s} \\[6pt]
% \dfrac{1}{s} \\[4pt]
% \vdots \\[4pt]
% \dfrac{1}{s}
% \end{pmatrix}.
% \]
The optimal SSP coefficient of an implicit method of order two is exactly $2s$, 
and it is achieved by a DIRK method: full implicitness gives no advantage.
Note that this is simply the implicit midpoint rule repeated 
$s$ times, with a scaled step size.
\end{example}

\section{A dictionary between SSP conditions and absorbing Markov chains}
\label{sec:dictionary}

The main conceptual contribution of AI to this proof was the use of 
absorbing Markov chains to prove the critical inequality in Lemma \ref{lem:AS}. 
ChatGPT introduced us to the connection to the matrices in SSP theory, 
which rest on the following translation.

% which was the main
% conceptual contribution of the AI tools to this work: the identification of the
% matrix $Q$ below as a \emph{substochastic} matrix was not known to us before.

Fix $r > 0$ satisfying \eqref{eq:KFSconditions} and define the following matrices 
as in \eqref{eq:NQw}
\begin{equation*}
  N = I + rA, \qquad Q = I - N^{-1} = rA(I+rA)^{-1}, \qquad w^\top  = r\, b^\top  N^{-1}.
\end{equation*}

Conditions \eqref{NSSP1}--\eqref{NSSP5}  imply
\begin{equation} \label{eq:Qcond}
  Q \ge 0, \qquad Q\one \le \one, \qquad w \ge 0, \qquad w^\top \one \le 1 .
\end{equation}
Moreover, $I-Q=N^{-1}$ is invertible, 
so Lemma \ref{lem:substochastic} gives  $\rho(Q)<1$. 
Thus $Q$ is the transition matrix among the transient states of an absorbing discrete-time Markov chain and $w$ is a (possibly defective) initial distribution on those states.  Then
$N = (I-Q)^{-1}$ is the fundamental matrix of the chain \cite{KemenySnell1976}, and
the absorption time $X$ of the associated continuous-time chain with transient generator $-(I-Q) = -N^{-1}$ and initial distribution $w$ 
is a \emph{phase-type}
random variable of order at most $s$ \cite{Neuts1981,Asmussen2003}, 
with an atom
$\Prob(X=0) = 1 - w^\top \one$ and moments
\begin{equation}\label{eq:PHmoments}
\E [X^0] = 1, \qquad  \E [X^k] \;=\; k!\; w^\top  N^k \one, \qquad k = 1,2,\dots
\end{equation}
%\lajos{for consistent notation throughout the manuscript I'd suggest here and below: $\E[X^k]$, $\E[X]$, $\E[X^2]$, $\E[X^3]$ , please check}
% \sigal{Zack and I did this, hopefully everywhere.}
Finally, the ``tall-tree'' order conditions translate into \emph{moment conditions}:
if the method has order at least $p$, then $b^\top A^{k-1}\one = 1/k!$ for
$1 \le k \le p$, and therefore, expanding $N^{k-1} = (I+rA)^{k-1}$,
\begin{equation}\label{eq:momentpin}
  w^\top  N^k \one
  \;=\; r\, b^\top  (I+rA)^{k-1}\one
  \;=\; r \sum_{j=0}^{k-1} \binom{k-1}{j} \frac{r^{\,j}}{(j+1)!},
  \qquad 1 \le k \le p .
\end{equation}
In particular,
\begin{equation}\label{eq:moments123}
  \E [X] = r, \qquad
  \E [X^2] = r^2 + 2r \quad (p \ge 2), \qquad
  \E [X^3] = r^3 + 6r^2 + 6r \quad (p \ge 3).
\end{equation}

\begin{center}
\begin{tabular}{ll}
\toprule
SSP / Runge--Kutta object & probabilistic object \\
\midrule
conditions \eqref{NSSP1}--\eqref{NSSP5} at $r$ & absorbing-chain representation $(Q,w)$ \\
$N = I + rA$ & fundamental matrix $(I-Q)^{-1}$ \\
SSP value $r$ & mean absorption time $\E [X]$ \\
order-$k$ tall-tree condition & $k^\text{th}$ moment condition \eqref{eq:momentpin} \\
optimal second-order method & Erlang distribution (serial chain) \\
\bottomrule
\end{tabular}
\end{center}

This translation was critically needed in the proof of the
substochastic sequence inequality in Lemma \ref{lem:AS};
the remainder of the Lemmas in Section \ref{sec:SSPbound} do not require the language  of absorbing Markov chains. However, it is possible to formulate the  entire proof in this language, which we do in the next section.

\subsection{A self-contained probabilistic proof}\label{sec:altproof}

The following alternative proof works directly with a discrete chain and
uniformization; it is somewhat cleaner in that Lemma \ref{lem:AS} and
Lemma \ref{lem:wcs} are absorbed into standard probability. 
We re-state the theorem with only the minimal conditions on the coefficients:

\begin{theorem}\label{thm:alt} (Re-statement of Theorem \ref{thm:SSPbound})
Let $N$ be an invertible real $s\times s$ matrix with $M = N^{-1}$ satisfying,
entrywise, $M\one\geq 0$ and $M\leq I$.  Let $w \ge 0$ with $w^\top \one\leq 1$, and
let $r\geq0$ satisfy $w^\top N\one=r$ and $w^\top N^2\one=r+\tfrac{r^2}{2}$.
Then $r\leq 2s$.
\end{theorem}

\begin{proof}
Set $Q=I-M$.  Then $Q\geq0$ and $Q\one\le \one$ (both understood componentwise), 
so $Q$ is substochastic, and since
$I-Q=M$ is invertible, Perron--Frobenius theory implies $\rho(Q)<1$ (as in Lemma
\ref{lem:substochastic}).  Consequently $N=(I-Q)^{-1}=\sum_{k\ge0} Q^k$ is the
fundamental matrix of the absorbing Markov chain with transient transition matrix
$Q$ \cite[Ch.~3]{KemenySnell1976}.

Start this chain in transient state $i$ with probability $w_i$, and in the absorbing
state with the remaining probability $1-w^\top \one$.  Let $T\in\{0,1,2,\ldots\}$ be
its absorption time.  For every $k\geq0$, $\Prob(T>k)=w^\top Q^k\one$.  The tail-sum
identity $n=\sum_{k\ge0} \mathbf 1_{\{n>k\}}$ and $n^2=\sum_{k=0}^{n-1}(2k+1)$ 
give, by Tonelli's theorem,
\[
    \E[T]=\sum_{k\ge0} \Prob(T>k)=w^\top N\one=r, \qquad
    \E[T^2]=\sum_{k\ge0}(2k+1)\Prob(T>k)=w^\top(2N^2-N)\one = r+r^2,
\]
where we used $\sum_k k Q^k = Q(I-Q)^{-2} = Q N^2 = N^2-N$.  Hence
\begin{equation}\label{eq:varT}
    \Var(T)=\E[T^2]-(\E[T])^2=r .
\end{equation}
%\lajos{for consistent notation throughout the manuscript I'd suggest here and below: $\E[T]$, $\E[T^2]$, please check}
%\sigal{Zack and I think we did this now}
These are the standard moment formulas of a discrete phase-type distribution
\cite[Ch.~2]{LatoucheRamaswami1999}.

We pass to continuous time by uniformization: let $E_1,E_2,\ldots$ be i.i.d.\ rate-one
exponentials independent of the chain and set $X=\sum_{j=1}^{T}E_j$.  Then $X$ is
the absorption time of the continuous-time chain with transient generator
$Q-I=-M$, hence a continuous phase-type variable with at most $s$ phases
\cite[Ch.~III]{Asmussen2003}.  Conditionally on $T$, $X$ is Erlang with $T$ rate-one
stages, so $\E[X\mid T]=T$ and $\Var(X\mid T)=T$, and by \eqref{eq:varT} and the
laws of total expectation and variance,
\[
    \E [X] = r, \qquad \Var(X) = \E [T]+\Var(T)=2r .
\]
Aldous and Shepp \cite{AldousShepp1987} proved that the Erlang distribution is the
least variable phase-type distribution of a given order: every continuous phase-type
$X$ with at most $s$ phases and $\E [X]>0$ satisfies $\Var(X)/(\E [X])^2\geq 1/s$.
If $r=0$ we are done; if $r>0$, then $2r \geq r^2/s$, i.e.\ $r\leq2s$.
\end{proof}

\begin{remark}[uniqueness of the optimal method]\label{rem:uniqueness}
Equality in the Aldous--Shepp inequality holds only for the Erlang distribution,
i.e.\ when the chain $(Q,w)$ is a serial chain traversed from its head.  
We can show that the $r = 2s$  is attained \emph{if and only if} 
the method is a permutation of the DIRK method above. 
The proof in Lean is given by
\leanproof{orderTwo\_two\_mul\_iff\_halfRateSerialMat}{SSPBoundEquality}{403}.

\end{remark}

\section{Extending the $2s$ Theorem using the connections between SSP methods and 
phase-type distributions}
%conjecture A moment hierarchy: towards sharp order-dependent bounds}
\label{sec:hierarchy}

The dictionary of Section \ref{sec:dictionary} converts 
the SSP conditions \eqref{eq:KFSconditions} at $r$ combined with
order-$p$ accuracy into a \emph{moment problem}: 
the phase-type absorption time $X$
of order at most $s$ must satisfy \eqref{eq:momentpin}.  
The bound $r\le 2s$ used a single moment inequality;  
higher moments may give more. In  Section \ref{sec:moment-hierarchy} 
we provide some results 
on a moment-ratio hierarchy for acyclic and cyclic phase-type distributions.
In Section \ref{sec:SSPboundP3} we use the results of Section 
\ref{sec:moment-hierarchy} with the 
dictionary of  Section \ref{sec:dictionary} to establish a bound on
SSP methods that satisfy the third order tall tree conditions.

\subsection{Moment-ratio hierarchy results} \label{sec:moment-hierarchy} 

We begin by citing two theorems from the literature: 
% \lajos{we mentioned the Aldous--Shepp theorem earlier several times and it is presented only here? a possible resolution: we could give some forward references to point to this section} \sigal{Done!}

\begin{theorem}[Aldous--Shepp least variability theorem \cite{AldousShepp1987}]
\label{thm:AldousShepp}
Let $s \ge 1$ be an integer, and let $X$ admit a
phase-type representation with at most $s$ transient states,
possibly with an atom at $0$. Let
\[m_0=1 , \qquad m_1 = {\E}[X] , \qquad m_2 = {\E}[X^2] ,\]
then 
\[ m_{0} m_{2} \ge \frac{s+1}{s} m_1^2. \]
\end{theorem}

\begin{theorem}[Osogami's least-$r$-value inequality for acyclic
phase-type distributions \cite{Osogami2005,OsogamiHarcholBalter2006}]
\label{thm:OsogamiHarcholBalter}
Let $s \ge 1$ be an integer, and let $X$ admit an {\em acyclic}
phase-type representation with at most $s$ transient states,
possibly with an atom at $0$. Write
$m_k={\E}[X^k]$, with $m_0=1$, and assume $m_1>0$.
Then
\[ \frac{m_{1}m_{3}}{m_2^2} \ge \frac{s+2}{s+1}. \]
\end{theorem}

% The case $k=2$ is the least-$r$-value inequality for acyclic
% phase-type distributions established in \cite[Lemma~3]{OsogamiHarcholBalter2006}
% following the thesis  \cite{Osogami2005}. 
We now extend Osogami's Theorem (Theorem \ref{thm:OsogamiHarcholBalter})
to  all $k \geq 1$:

\begin{theorem}[Moment-ratio hierarchy for acyclic phase-type distributions]
\label{thm:acyclic-hierarchy}
Let $s \ge 1$ be an integer, and let $X$ admit an acyclic
phase-type representation with at most $s$ transient states,
possibly with an atom at $0$. Write
$m_k={\E}[X^k]$, with $m_0=1$, and assume $m_1>0$.
Then, for every integer $k \ge 1$,
\[
\frac{m_{k-1}m_{k+1}}{m_k^2}
\ge \frac{s+k}{s+k-1}.
\]
\end{theorem}

\begin{proof}
Conditional on its path through the transient states, $X$ is
a sum of at most $s$ independent exponential holding times,
since acyclicity rules out repeat visits.
Consequently, we may represent 
%\lajos{if this $\overset{d}{=}$ is a defining equality, I suggest using $:=$}
\[
X {:=} \sum_{i=1}^s A_i E_i,
\]
where $E_1,\ldots,E_s$ are independent unit-rate exponential
random variables, and $A_i$ is the reciprocal exit rate of the
$i^\text{th}$ visited state. Because the unit exponentials driving the
holding times are independent of the embedded jump chain, the
nonnegative random vector $(A_1,\ldots,A_s)$ is independent of
them. Unused stages have coefficient zero, and immediate absorption
corresponds to all coefficients being zero (we use the convention that 
$0^0=1$).
Note that $A_i\le 1/\lambda_{\min}$, where $\lambda_{\min}>0$ is
the smallest exit rate.
Set $G=\sum_{i=1}^s E_i$ to define the normalization
\[ U_i=\frac{E_i}{G}.\]
The gamma--Dirichlet decomposition shows that
$G\sim {\Gamma}(s,1)$ is independent of
$(U_1,\ldots,U_s)$. Since $(G,U)$ is independent of $A$, it follows
that $G$ is independent of $(A,U)$. We can therefore express $X$ as
\[
X {:=} GZ,
\qquad Z=\sum_{i=1}^s A_iU_i,
\]
with $G$ and $Z$ independent, and with $0\le Z\le 1/\lambda_{\min}$,
so that all moments are finite. Writing
\[ (s)_j:=\frac{\Gamma(s+j)}{\Gamma(s)}
=s(s+1)\cdots(s+j-1),\qquad (s)_0=1, \] we obtain
\[
m_j= {\E}[X^j] = (s)_j\,\mathbb{E}[Z^j].
\]
These quantities are positive for all $j$, since $m_1>0$ forces
$\mathbb{P}(X>0)>0$.
By Cauchy--Schwarz we have 
\[ \left({\E}[Z^k]\right)^2 \leq {\E}[Z^{k-1}] {\E}[Z^{k+1}], \]
so that 
\[  \frac{m_{k-1}m_{k+1}}{m_k^2} =
\frac{s+k}{s+k-1}
\frac{{\E}[Z^{k-1}]{\E}[Z^{k+1}]}{\left({\E}[Z^k]\right)^2}
\ge \frac{s+k}{s+k-1}.
\]
\end{proof}

%\noindent \sigal{Barry: let's think about whether we need the equality result. For now
%it is commented out below.}
% \begin{remark}[The equality case]
% For $k=1$, equality holds if and only if $X$ has an Erlang
% distribution with $s$ stages. For any fixed $k \ge 2$, equality
% holds if and only if 
% \[ \mathcal{L}(X) =(1-p)\delta_0+p\,\operatorname{Erlang}(s,\lambda)
% \] for some $p\in(0,1]$ and $\lambda>0$.
% To see this, notice that 
% For $k=1$, equality in Cauchy--Schwarz requires $Z$ to be
%  a positive constant almost surely. For $k\ge2$, it requires
%  $Z$ to be a positive constant conditional on $Z>0$.
%  These conditions yield the stated Erlang distributions
%  and their mixtures with an atom at zero. Conversely,
%  direct calculation verifies equality for those distributions.
% \end{remark}

\begin{remark}[Results for the cyclic case]
As we mentioned above, for the 
the case $k=1$ this result is the Aldous--Shepp
least-variability Theorem \ref{thm:AldousShepp} 
which holds for the cyclic as well as acyclic 
cases \cite{AldousShepp1987}. 
For $k=1$ we do not require that  $X$ admit an {\em acyclic}
phase-type representation with at most $s$ transient states.
However, Osogami's Theorem \ref{thm:OsogamiHarcholBalter} which gives the 
result for $k=2$ in \cite{OsogamiHarcholBalter2006} does require
an acyclic phase-type distribution, and our extension to all $k\geq 1$ in
Theorem \ref{thm:acyclic-hierarchy} does as well.
We are able to extend this result to the {\em cyclic case} for $k=2$, 
which will be useful to us in the next section, where we prove a bound on the SSP
coefficient for third order methods. Recall that 
in our SSP language {\em acyclic} is equivalent to $A$ being a triangular matrix
after simultaneous permutation, while a truly full matrix $A$ that corresponds 
to the {\em cyclic} phase-type representation.
\end{remark}

\begin{theorem}[Moment-ratio hierarchy for cyclic phase-type distributions]
\label{thm:cyclic-hierarchy}
Let $s \ge 1$ be an integer, and let $X$ admit a
phase-type representation with at most $s$ transient states,
possibly with an atom at $0$. Write
$m_k={\E}[X^k]$, with $m_0=1$, and assume $m_1>0$.
Then,  we have
\begin{eqnarray} \label{eq:moment}
\frac{m_{1}m_{3}}{m_2^2} \ge \frac{s+2}{s+1}.
\end{eqnarray} 
\end{theorem}
\noindent Before we prove this theorem, we state the following Lemma, the proof
of which is algebraic and given in Appendix \ref{app:LemAlg}. 
\begin{lemma}[Algebraic inequality]\label{lem:alg}
Let $n\ge 2$, $w\in[0,1)$, $u\ge0$, and
$\beta\ge\frac{n}{n-1}\,w^2$; we define $\beta^2=0$ when 
$w=0$, and in particular $\beta^2/w = 0$  for $w=0$. 
% \lajos{I don't understand: with $\beta=0$ if $w=0$, and a verb is missing from the second half of the sentence? I think we want to define $\beta^2/w=0$ for $w=0$?} 
Put
\[ m = 2(1-w)+u+\beta . \]
Then
\begin{equation}\label{eq:alg}
  3(1-w)\,m + 3u +  \frac{n+1}{n} \left(\frac{u^2}{1-w}+\frac{\beta^2}{w}\right)
  \;\ge\; \frac{n+2}{n+1}\,m^2 .
\end{equation}
\end{lemma}

\begin{proof}[Proof of Theorem \ref{thm:cyclic-hierarchy}.]
If \eqref{eq:moment} holds for every fixed starting state, then by Cauchy–Schwarz it also holds for any random initial state; thus it is enough to prove the theorem for a chain starting at a fixed transient state.
% Write $m_k=\sum_i \alpha_i m_k^{(i)}$ for $k\ge 1$; the atom at $0$ contributes nothing. With $c=\frac{s+2}{s+1}$ and $m_2^{(i)}\le c^{-1/2}\sqrt{m_1^{(i)}m_3^{(i)}}$ for each $i$, Cauchy–Schwarz gives
% $$
% m_2=\sum_i \alpha_i m_2^{(i)}\le c^{-1/2}\sum_i \sqrt{\alpha_i m_1^{(i)}}\,\sqrt{\alpha_i m_3^{(i)}}\le c^{-1/2}\sqrt{m_1 m_3},
% $$
% and squaring yields \eqref{eq:moment}.
The proof is based on an induction argument:

For the base case $n=1$ (one transient state) the absorption time is 
exponential with some rate \(\lambda>0\), so
\[
m_1=\frac1\lambda,\qquad
m_2=\frac2{\lambda^2},\qquad
m_3=\frac6{\lambda^3},
\qquad
\Rightarrow
\qquad
\frac{m_1m_3}{m_2^2} =\frac64=\frac32 =\frac{n+2}{n+1}. \] 
The equality still holds if there is an atom at zero.

Now we continue the induction argument:
Let $n\ge2$ and assume that the identity 
\[  n \;  m_1 m_3 \ge (n+1) m_2^2 \]
holds for the moments $m_j$ of any phase-type distribution with at most $n-1$ states.

\noindent {\bf Step 1: Decomposing the path.}
Take an $n$-state chain started at a fixed transient state $i$, and let
$Z$ be its absorption time. Decompose the path at the first departure from $i$:
\begin{itemize}
  \item $E$ is the initial holding time at $i$; it is exponential with
        mean $x>0$.
  \item After leaving $i$, the chain returns to $i$ before absorption
        with probability $p$, and is absorbed first with probability
        $q=1-p$. Since absorption is almost sure, $q>0$.
  \item $U$ is the time from departure until return to $i$,
        conditioned on returning (only needed when $p>0$).
  \item $V$ is the time from departure until absorption, conditioned on
        no return. ($V$ may have an atom at $0$, if the chain jumps from
        $i$ straight to absorption.)
\end{itemize}
Run the chain from the
moment it leaves $i$, and stop it when it first hits either $i$ or the
absorbing state. Only the $n-1$ states other than $i$ are visited in
between. 
 Conditioning on which target is hit first (a Doob
 $h$-transform) changes the transition rates among these $n-1$ states but
 adds no new ones. 
\emph{So both $U$ and $V$ are phase-type with at most $n-1$ transient states.} 

\noindent {\bf Step 2: Renewal equation.}
For a random variable $Y$ write $F_Y(t)=\E[e^{tY}]$; 
this is finite for $t$ near $0$ since all variables here
are phase-type. 
\begin{itemize}
\item {\bf Holding time:} The chain first stays at $i$ for a time $E$,
        exponential with mean $x$, so \[F_E(t)=1/(1-xt).\] In a
        continuous-time Markov chain, $E$ is independent of everything
        that happens after the departure.
\item {\bf Return} (probability $p$): 
In this case the excursion returns to $i$ after time $U$. 
By the strong Markov property the process then restarts, 
and the remaining time $Z'$ to absorption is a copy of
$Z$, independent of $E$ and $U$. 
So in this event we have
\[ Z = E + U + Z' \; \; \;  \mbox{with probability} \; \; \; p. \] 
\item  {\bf No return} (probability $q$): If the chain does not return, 
the excursion ends in absorption after time $V$, and we have 
$E$, then $V$, followed by absorption
\[ Z = E+V \qquad  \mbox{with probability} \; \; \;  q. \]
\end{itemize}
Taking expectations of $e^{tZ}$ over the two events, and using the fact that
moment generating functions  of independent sums multiply, we get
\[  F_Z(t) = p\,F_E(t)F_U(t)F_Z(t) + q\,F_E(t)F_V(t)
         = \frac{p\,F_U(t)F_Z(t)+q\,F_V(t)}{1-xt}, \]
or equivalently
\[ \bigl(1-xt-p\,F_U(t)\bigr)\,F_Z(t) \;=\; q\,F_V(t). \]
When $p=0$, drop every term involving $U$.

\noindent {\bf Step 3: Normalize.}
Since $\E[Z]\ge x>0$ we can normalize $Z$ so that $\E[Z] = 1$. This does not change the 
identity because  $\lambda Z$ ($\lambda>0$) multiplies $\E[Z]\,\E[Z^3]$ and $\E[Z^2]^2$
by the same factor $\lambda^4$. Rescaling corresponds to multiplying all
rates of the chain by $1/\lambda$. This normalization
rescales $E$, $U$, $V$ by the same factor, and leaves
$p$ and $q$ unchanged. 

\noindent {\bf Step 4: Moment identities.} Divide 
\[ \bigl(1-xt-p\,F_U(t)\bigr)\,F_Z(t) \;=\; q\,F_V(t). \]
by $q$ and expand in powers of $t$. Using
$(1-p)/q=1$ and $F_U(t)=1+t\,\E[U]+\tfrac{t^2}{2}\E[U^2]
+\tfrac{t^3}{6}\E[U^3]+O(t^4)$, we have 
\[
  \frac{1-xt-p\,F_U(t)}{q}
  = 1 - \frac{x+p\,\E[U]}{q}\,t
      - \frac{t^2}{2}\,\frac pq\,\E[U^2]
      - \frac{t^3}{6}\,\frac pq\,\E[U^3]
    + O(t^4).
\]
Multiply this by $F_Z(t)=1+t+\tfrac{t^2}{2}\E[Z^2]+\tfrac{t^3}{6}\E[Z^3]
+O(t^4)$, and match the coefficients of $t$, $t^2$ and $t^3$ with those
of $F_V(t)=1+t\,\E[V]+\tfrac{t^2}{2}\E[V^2]+\tfrac{t^3}{6}\E[V^3]+O(t^4)$.
The coefficient of $t$ gives
\begin{equation}\label{eq:first}
  \frac{x+p\,\E[U]}{q} \;=\; 1-\E[V].
\end{equation}
Since $x>0$, this equation gives $1-\E[V]>0$, so
$0\le \E[V]<1$. Now the coefficients of $t^2$ and $t^3$
give
\begin{subequations}
\begin{align}
  \E[Z^2] &= 2\bigl(1-\E[V]\bigr) + \frac pq\,\E[U^2] + \E[V^2],
  \label{eq:second}\\
  \E[Z^3] &= 3\bigl(1-\E[V]\bigr)\E[Z^2] + 3\,\frac pq\,\E[U^2]
            + \frac pq\,\E[U^3] + \E[V^3].
  \label{eq:third}
\end{align}
\end{subequations}

\noindent\textbf{Step 5: Apply Lemma~\ref{lem:alg}.}
Since $U,V$ are phase-type with at most $n-1$ states, 
they satisfy
\[
  \E[V] \E[V^3]\;\ge\; \frac{n+1}{n} \left(\E[V^2]\right)^2, \qquad
   \E[U] \E[U^3]\;\ge\; \frac{n+1}{n} \left(\E[U^2]\right)^2
\]
First suppose  $p>0$,  then multiplying the second inequality by
$p/q$ and dividing by $\E[U]$ gives
\[   \frac pq\,\E[U^3]
  \;\ge\; \frac{n+1}{n}\,\frac{\bigl(\frac pq\E[U^2]\bigr)^2}{\frac pq\E[U]}
  \;\ge\;  \frac{n+1}{n}\,\frac{\bigl(\frac pq\E[U^2]\bigr)^2}{1-\E[V]} ,
\]
where we used  $\frac pq\E[U] = 1-\E[V]-\frac xq < 1-\E[V]$, which
follows from \eqref{eq:first} and $x>0$. If $p=0$, both sides vanish,
so this inequality still holds.
We now return to \eqref{eq:third} 
\begin{eqnarray*}
\E[Z^3] & = & 3\bigl(1-\E[V]\bigr)\E[Z^2] + 3\,\frac pq\,\E[U^2]
            + \frac pq\,\E[U^3] + \E[V^3]  \\
& \geq  &   3\bigl(1-\E[V]\bigr)\E[Z^2] + 3\,\frac pq\,\E[U^2]
         + \frac{n+1}{n}\,\frac{\bigl(\frac pq\E[U^2]\bigr)^2}{1-\E[V]}  
         + \E[V^3] \\
  & \geq  &      3\bigl(1-\E[V]\bigr)\E[Z^2] + 3\,\frac pq\,\E[U^2]
  +  \frac{n+1}{n} \left(
      \frac{\bigl(\frac pq\E[U^2]\bigr)^2}{1-\E[V]}
      +\frac{(\E[V^2])^2}{\E[V]}
    \right).    
\end{eqnarray*} 
Now if we define 
\[   w=\E[V],\qquad u=\frac pq\,\E[U^2],\qquad \beta=\E[V^2], \]
then  the hypothesis of Lemma~\ref{lem:alg} hold
\[ 0\le w<1, \qquad  u\ge0, \qquad
\beta\ge\frac{n}{n-1}w^2.\]
The final inequality is simply  Aldous--Shepp  for $V$, and $\beta=0$ when $w=0$.
The value  $m=2(1-w)+u+\beta = \E[Z^2]$ in Lemma~\ref{lem:alg}
gives 
\[   \E[Z^3]\;\ge\;\frac{n+2}{n+1}\,\left(\E[Z^2]\right)^2 ,\]
which, since $\E[Z] = 1$ is 
\[  \E[Z] \E[Z^3]\;\ge\;\frac{n+2}{n+1}\,\left(\E[Z^2]\right)^2 .\]
\end{proof}

Given Theorem \ref{thm:acyclic-hierarchy} that establishes the relationship  
\[
\frac{m_{k-1}m_{k+1}}{m_k^2}
\ge \frac{s+k}{s+k-1}
\]
in the acyclic case for every integer $k \ge 1$ and $s\geq 1$, and the promising
result for the cyclic case in Theorem \ref{thm:cyclic-hierarchy} 
for integers $1 \leq k \leq 2$ and  $s\geq 1$ it is tempting to conjecture that this 
relationship is true for all $k \geq 1$ and $s\geq 1$ in the cyclic case:

\begin{conjecture} \label{conj:cyclic-hierarchy}
Let $s \ge 1$ be an integer, and let $X$ admit a cyclic
phase-type representation with at most $s$ transient states,
possibly with an atom at $0$. Write
$m_k={\E}[X^k]$, with $m_0=1$, and assume $m_1>0$.
Then, for every integer $k \ge 1$,
\[
\frac{m_{k-1}m_{k+1}}{m_k^2}
\ge \frac{s+k}{s+k-1}.
\]
\end{conjecture}
\noindent  However, this conjecture is false. 
Consider the counterexample 
for $s=6, k=71$:
\begin{example}
Let $\varepsilon = 10^{-8}$ and
\[ Q =
\begin{pmatrix}
\tfrac19 & \tfrac89 & 0 & 0 & 0 & 0\\[2pt]
0 & \tfrac19 & \tfrac89 & 0 & 0 & 0\\[2pt]
0 & 0 & \tfrac19 & \tfrac89 & 0 & 0\\[2pt]
0 & 0 & 0 & \tfrac19 & \tfrac89 & 0\\[2pt]
0 & 0 & 0 & 0 & \tfrac19 & \tfrac89\\[2pt]
\varepsilon & 0 & 0 & 0 & 0 & 0
\end{pmatrix}.
\]
Consider the phase-type absorption time with transient generator $Q - I$,
started in state~$6$. This is admissible: $Q \ge 0$ entrywise and
$Q\one \le \one$. 
Moreover, $I-Q$ is a nonsingular $M$-matrix, with
\[ \det(I-Q) = \frac{1422222208}{2562890625} \neq 0
\qquad\text{and}\qquad
N := (I-Q)^{-1} \ge 0,
\]
so absorption occurs almost surely. The only exit rate is from state~$6$,
namely $\bigl((I-Q)\one\bigr)_6 = 1-\varepsilon$.

The moments are
\[ m_j = j!\, e_6^{\top} N^{j}\one, \qquad j \ge 1 . \]
At $k = 71$, if the conjecture  were true, we would have 
\[ \frac{m_{70}\, m_{72}}{m_{71}^{2}} \geq \frac{6 + 71}{6 + 71 - 1} = \frac{77}{76}.\]
However, we can compute
\[ \frac{m_{70}\, m_{72}}{m_{71}^{2}}
= 1.013156577010572151\ldots
\;<\;
\frac{77}{76} = 1.013157894736842105\ldots,
\]
with gap
\[
\frac{m_{70}\, m_{72}}{m_{71}^{2}} - \frac{77}{76}
\approx -1.31772627\times 10^{-6} < 0.
\]
The sign of the gap is certified by exact rational arithmetic.
\end{example}

% \sigal{Shouold we put the proof in? Refer to it and write it elsewherew? I have not beem
% able to independently understand either of the proofs -- Barry and I will discuss.
% }
% \revnote{Literature check needed: the $k = 1$ case is Aldous--Shepp; is the general
% case known?  Candidate sources: O'Cinneide's phase-type majorization work
% \cite{OCinneide1991}, and the moment-matching literature for acyclic phase-type
% distributions \cite{BobbioHorvathTelek2005}, where the order-$s$ moment regions are
% characterized for APH.}
\begin{remark}
    In the process of getting these results for the cyclic case, 
    we generated several other proofs that may be interesting in the 
    area of probability but are not directly related to the SSP Runge--Kutta methods
    results we are interested in. We provide the writeup and Lean verification at \leanreporef, and list these results here. The results listed below
    are all AI-derived  and Lean verified, however they were not all verified by humans.
   % \sigal{BARRY! Help here please:}
    \begin{enumerate}
        \item We extended Theorem \ref{thm:cyclic-hierarchy} to the case $k=3$  
        i.e. to the bound
     \[ \frac{m_{2}m_{4}}{m_3^2} \ge \frac{s+3}{s+2}. \] The proof is described in \leandocfile{k3\_all\_stages\_induction\_summary.pdf} and the final proof in Lean is \leanproof{hierarchy\_k3}{K3Complete}{40}.
     \item We are able to extend the results in Theorem \ref{thm:cyclic-hierarchy}
     for all $k$ so that 
     \[  \frac{m_{k-1}m_{k+1}}{m_k^2} 
\ge \frac{s+k}{s+k-1} \mbox{ for } s=2,3 .\]
    The proofs are described in \leandocfile{cyclic\_all\_k\_fisher\_reduction.pdf} and the final theorems in Lean are \leanproof{momentConjecture\_two}{MatrixMomentReduction}{370} for $s=2$ and
    \leanproof{conjecture\_one\_three\_state} {ThreeStateConjectures}{109} for $s=3$.
     \item  We can further extend the results in Theorem \ref{thm:cyclic-hierarchy}
     for the case where $A$ (and therefore $N$) are block triangular matrices with 
     blocks of size $s \leq 3$.  This follows from the previous item and the result that conjecture
     \ref{conj:cyclic-hierarchy} holds for a block triangular matrix if it holds for each block on the diagonal.
     The proofs are described in 
     \leandocfile{cyclic\_all\_k\_fisher\_reduction.pdf} 
     and the final theorem in Lean is 
     \leanproof{conjecture\_one\_block\_triangular}{BlockTriangularMoments}{67}.
    \end{enumerate}
\end{remark}

\subsection{A bound on the SSP coefficient for third-order methods} \label{sec:SSPboundP3}
In this subsection we use the theorems above with $k=2$ to get a bound on the 
SSP coefficients of methods that satisfy the third order tall tree conditions.

\begin{theorem} \label{thm:SSPboundP3}
Let a Runge-Kutta method $(A,b)$ be given and define
matrices $N \in \mathbb{R}^{s\times s}$ and  $w \in \mathbb{R}^{s}$ 
by \[ N= I + rA, \; \; \;  w^\top  = r b^\top  N^{-1}.\] 
Suppose that they satisfy the SSP conditions \eqref{NemethSSPconditions}:
\[    N \; \; \; \mbox{is invertible} , \qquad
    N^{-1} \one \geq  0 ,  \qquad
    Q = I -  N^{-1}  \geq 0  , \]
\[    w^\top  \geq 0 , \qquad   1- w^\top  \one \geq 0, \]
(where the inequalities are understood componentwise),
and the order conditions 
\[w^\top  N \one = r , \qquad   w^\top  N^2 \one  =  r + \frac{1}{2} r^2, 
\qquad   w^\top  N^3 \one  = r + r^2 + \frac{1}{6} r^3. \]
 Then 
 \begin{equation} \label{eq:3rdObound}
     r \;\le\; s - 1 + \sqrt{s^2-1} .
 \end{equation}  
\end{theorem}

\begin{proof}
First, we note that for $r \leq 0$ the bound \eqref{eq:3rdObound}
is trivially satisfied for all 
$s \geq 1$. 
Therefore, suppose $r > 0$. If $X$ is the absorption time of the chain $(Q,w)$
where $Q = I - N^{-1}$ then it is a  phase-type of order at most $s$ with
moments \eqref{eq:PHmoments}
\[m_1 = r,
\qquad
m_2 = 2\,w^\top N^2\one = r^2+2r,
\qquad 
m_3 = 6\,w^\top N^3\one = r^3+6r^2+6r.\]
By  Theorem \ref{thm:acyclic-hierarchy} for the case where $A$
is triangular, or Theorem  \ref{thm:cyclic-hierarchy}  for the case where $A$ is not triangular,
we have
\[(s+1)\, m_1 m_3 \ge (s+2)\, m_2^2\]
which becomes
\[
  (s+1)\,r\,(r^3+6r^2+6r) \;\ge\; (s+2)\,(r^2+2r)^2 \quad
  \Rightarrow \quad
  r^2\bigl(r^2 - 2(s-1)r - 2(s-1)\bigr) \le 0.
\]
Dividing by
$r^2>0$ and completing the square, we get
\[(r-(s-1))^2 \le s^2-1 \qquad \Rightarrow \qquad
r \;\le\; s - 1 + \sqrt{s^2-1} . \]
\end{proof}

\noindent This theorem shows that every irreducible 
$s$ stage implicit Runge--Kutta method of 
order $p \ge 3$ has \[ \sspc(A,b) \le s-1+\sqrt{s^2-1}. \]  
In fact, this bound is sharp: for any $s\ge 2$,  
the $s$-stage SDIRK method with coefficients
\[ b_i=\frac1s,\qquad
a_{ij}=
\begin{cases}
\displaystyle \frac{1}{\sqrt{s^2-1}}, & j<i,\\[6pt]
\displaystyle \frac12\left(1-\sqrt{\frac{s-1}{s+1}}\right), & j=i,\\[6pt]
0, & j>i
\end{cases}
\]
has SSP coefficient $\sspc(A,b) = s-1+\sqrt{s^2-1}$
satisfies {\em all} the third order conditions:
\[
b^\top\mathbf{1}=1,\qquad
b^\top c=\frac12,\qquad
b^\top(c\circ c)=\frac13,\qquad
b^\top Ac=\frac16.
\]

In fact, the bound in Theorem \ref{thm:SSPboundP3} 
does not require all the order conditions, 
only the tall tree conditions, so that this is true even if the 
Runge--Kutta method only has  linear order $p_{lin} = 3$.
Methods that achieve this bound were found in \cite{FerracinaSpijker2008}
in an optimization over singly diagonally implicit Runge--Kutta (SDIRK)
methods. By definition, SDIRK methods have a lower triangular $A$ and 
all the elements on the diagonal of $A$ are the same.
In that work it was conjectured 
(\cite[Conjecture 3.2]{FerracinaSpijker2008})
that these methods gave the  optimal SSP coefficient
among SDIRK methods.

In \cite{KMG09}, these methods up to $s=11$ were found using optimization, 
this time in a search among {\em all} implicit Runge--Kutta methods.
The two-stage third order method was proven optimal by the fact that it
achieved the  threshold factor reported in \cite{vandegriend1986}.
The three-stage third order method was optimal in that
BARON  provided a certificate of global optimality \cite{KMG09}.  
In \cite{KMG09} the conjecture in \cite{FerracinaSpijker2008}
was extended: it was conjectured 
that the optimal  implicit SSP Runge–Kutta methods with 
$s\geq 2$ stages are  SDIRK,  and identical to these methods
\cite[Conjecture 2]{KMG09}.
This is strong evidence that the $k=2$ inequality is \emph{tight} 
for third order Runge--Kutta methods.
% Numerical evidence indicates that the $k=2$ inequality is \emph{tight} for third order Runge--Kutta methods; direct maximization of $ r$ under the third order tall tree conditions
% returns $r_2(s)$ to eight digits.

However, if we look at the positive root $r$ from the 
relationship between the $k=2,3,4$ moments for the acyclic case:
\[ (s+2)\,m_{2}(r)\,m_{4}(r) = (s+3)\,m_3(r)^2 \]
and compare it to  the numerically optimized
DIRK methods that satisfy the tall tree conditions to order $4$ we see that the 
for fourth order Runge--Kutta methods the  $k=3$ inequality is slack. 
This suggests that for higher orders the joint phase-type moment body
rather than the consecutive-triple log-convexity may be the binding object.
%(see \texttt{data/exp5\_summary.json}). The results above that show that
%$\sspc_{s,p_{lin} = 4} < r_3(s)$ suggest that this is indeed the case.

\section{An order barrier for algebraic stability with non-negative coefficients}
\label{sec:barrier}

Strong stability preservation is an extraordinarily strong property, and is particularly useful when dealing with the maximum norm, total variation semi-norm, or positivity.  For many problems, the essential bound involves instead an inner-product norm, which may lead
to looser conditions. Inner-product stability is often ensured
through a contractivity property:
\begin{align}\label{ipcontractivity}
    \langle f(u) - f(v), u-v\rangle \le 0.
\end{align}
A method is said to be B-stable if it guarantees the discrete contractivity property
\begin{align}
    \|u_{n+1} - v_{n+1}\| \le \|u_n - v_n\|
\end{align}
for all step sizes $\dt\ge 0$, for all contractive problems.
In order for a method to be B-stable, the following \emph{algebraic stability} condition is sufficient.  This condition is also necessary if the method has all abscissas distinct
\cite[Sec.~IV.12]{HairerWanner}.
\begin{definition}[\bf Algebraic Stability]
Let a Runge--Kutta method be given with coefficients $A, b$, and
let $B=\text{diag}(b)$. The method is algebraically stable if 
\begin{align}\label{algebraic-stability}
    B & \succeq 0 \qquad \qquad \text{ and } & M & 
    = BA + A^\top B - b b^\top \succeq 0.
\end{align}
(The notation $\succeq $ means the matrix is positive semidefinite.)
\end{definition}

This section was motivated by numerical searches for implicit Runge--Kutta methods that have non-negative coefficients (which are required for the SSP condition to
hold) and are algebraically stable.
When extensive numerical optimization failed to uncover any such methods of order higher than two, we began to wonder if the conditions for order three 
force a conflict between SSP and algebraic stability. This is explored in this section.

Once again, we are dealing with a Runge--Kutta method $(A,b)$.
The  order conditions for up to third order are:
\begin{subequations}  \label{eq:RKOC}
\begin{eqnarray}
b^\top \one    & = &  \sum_{i=1}^s b_i  =   1, \label{eq:RKOC1} \\
b^\top c & = & \sum_{i=1}^s  b_i c_i  =  \tfrac12,  \label{eq:RKOC2} \\
b^\top (c \circ c ) & = & \sum_{i=1}^s  b_i c_i^2  = \tfrac13, \label{eq:RKOC3a} \\
b^\top A c & = & \sum_{i=1}^s  \sum_{j=1}^s b_i a_{ij}   c_j  = \tfrac16. 
\label{eq:RKOC3b} 
\end{eqnarray}
\end{subequations}
In this section we show that a Runge--Kutta method with 
non-negative coefficients $(A, b)$ that satisfy the
order  conditions \eqref{eq:RKOC1}, \eqref{eq:RKOC2}, and \eqref{eq:RKOC3a} cannot also satisfy the  algebraic stability condition \eqref{algebraic-stability}.

\begin{theorem}[Order barrier for algebraic stability]\label{thm:barrier}
There is no Runge--Kutta method $(A,b)$ with, componentwise,
\[ A\ge 0,   \qquad   b \ge 0, \qquad A\one = c \]
satisfying the bushy tree conditions up to third order
\[  b^\top \one = 1, 
\qquad 
b^\top c = \tfrac12, \qquad
 b^\top (c \circ c)  =  \frac13 
\]
 whose algebraic-stability matrix
\[ M = B  A+  A^\top  B - b b^\top , \qquad \mbox{where} \; \; B= \diag(b), \]
is positive semidefinite.  
\end{theorem}

\begin{remark}\label{rem:cjnondecreasing}
For a Runge--Kutta method given by $A,b$, we can always
permute the rows and columns of  $A$ and the entries of  $b$
so that the abscissas $c_i$ are non-decreasing.
This type of permutation  merely relabels the stages and 
leaves the method unchanged in terms of its stability properties.
In our case we are not requiring that the method be diagonally
implicit  or explicit, in which case we would require the matrix $A$ to remain lower triangular or strictly lower triangular, so we
can freely perform this permutation.
Thus, in the following we  assume 
\[ c_1 \le c_2 \le \cdots \le c_s \]
without loss of generality.
\end{remark}

\noindent{\bf Outline of the proof:}
To prove this theorem, we will  show that if $M \succeq 0$ and $A, b \ge 0$
with $A\one =  c$, $b^\top \one = 1$, $b^\top c = \tfrac12$, then
$ b^\top (c \circ c)  < \frac13 $
This means that if the method has non-negative coefficients and
an algebraic stability matrix that is positive semidefinite then
the third order bushy tree condition cannot be satisfied.

\begin{lemma}\label{lem:partialsums}
Let $c_1\le c_2\le \cdots \le c_s$  and let
$b_i\ge 0$ satisfy
\[  \sum_{i=1}^s b_i=1,  \qquad  \sum_{i=1}^s b_i c_i=\frac12,
\qquad \sum_{i=1}^k b_i c_i  \ge 
\frac12 \left(\sum_{i=1}^k b_i \right)^2 \quad \mbox{for } \; k=1,\dots,s.\]
Then
\[  \sum_{i=1}^s b_i c_i^2  \le  \frac13-\frac1{12}\sum_{i=1}^s b_i^3   <  \frac13 .\]
\end{lemma}

\begin{proof}
Let
\[
\beta_0=0,\qquad \beta _k=\sum_{i=1}^k b_i,
\qquad m_i=\frac{1}{2} (\beta_{i-1}+\beta_i).
\]
Note that
\[ b_i m_i
= (\beta_i-\beta_{i-1})\left(\frac{\beta_i+\beta_{i-1}}{2}\right)
=\frac{1}{2} \left( \beta_i^2-\beta_{i-1}^2\right),
\]
so that
\[
\sum_{i=1}^k b_i m_i=
\frac{1}{2} \left( \beta_k^2-\beta_{k-1}^2\right)
+ \frac{1}{2} \left( \beta_{k-1}^2-\beta_{k-2}^2\right)
+ ... + \frac{1}{2} \left( \beta_{1}^2-\beta_{0}^2\right)
=\frac{1}{2} \beta_k^2.
\]
Now define $\delta_0 = 0$ and 
\[
\delta_k=\sum_{i=1}^k b_i(c_i-m_i), \qquad k = 1,..., s,
\]
so that
\[
\delta_k=\sum_{i=1}^k b_i c_i-  \sum_{i=1}^k b_i m_i
\geq  \frac12 \left(\sum_{i=1}^k b_i \right)^2 -  \frac{1}{2} \beta_k^2
\geq 0,\]
using the assumptions of the lemma.
The order conditions \eqref{eq:RKOC1} and \eqref{eq:RKOC2} 
(also called the {\em mean conditions}) give the final value
\[
\delta_s=\sum_{i=1}^s b_i c_i-\frac{\beta_s^2}{2}
=\frac12-\frac12=0.
\]

Moreover,
\[
m_{i+1}-m_i= \frac{1}{2} (\beta_{i}+\beta_{i+1}) - \frac{1}{2} (\beta_{i-1}+\beta_i) = \frac{1}{2} \left(\beta_{i+1} - \beta_{i-1} \right)
= \frac{b_i+b_{i+1}}2 \ge0,
\]
so the sequence $\left\{ m_i\right\}$ is non-decreasing.
Let  $ r_i=c_i+m_i ,$
so we can write 
\[
\sum_{i=1}^s b_i(c_i^2-m_i^2)
=\sum_{i=1}^s r_i\,b_i(c_i-m_i)
=\sum_{i=1}^s r_i(\delta_i-\delta_{i-1})
=\sum_{i=1}^{s-1} r_i \delta_i
  \;-\; \sum_{i=1}^{s-1} r_{i+1}\,\delta_i,
\]
using the fact that $\delta_0 = \delta_s  = 0$.
Thus we have
\[ \sum_{i=1}^s b_i(c_i^2-m_i^2) 
= - \sum_{i=1}^{s-1} \delta_i \left(r_{i+1} - r_i \right) \leq 0,\]
because $\delta_i \geq 0$ and  $r_i$ are non-decreasing (as it is the sum of 
non-decreasing  $c_i$  and non-decreasing $m_i$).
Consequently,
\[
\sum_{i=1}^s b_i c_i^2\le \sum_{i=1}^s b_i m_i^2.
\]

Finally, observe that   \(m_i\) are the midpoints of the consecutive intervals
\([\beta_{i-1},\beta_i]\)  of lengths \(b_i\), that  partition \([0,1]\).
We write
\begin{eqnarray*}
\int_{\beta_{i-1}}^{\beta_i}x^2\,dx
&=& \frac{1}{3}\left(\beta_i^3 - \beta_{i-1}^3 \right) 
= \frac{1}{3} (\beta_i- \beta_{i-1})( \beta_i^2 + \beta_i \beta_{i-1} + \beta_{i-1}^2 )\\
&=& \frac{1}{12} b_i ( 4 \beta_i^2 + 4 \beta_i \beta_{i-1} + 4 \beta_{i-1}^2 )\\
&=& \frac{1}{12} b_i
\left( 3 (\beta_i^2 + 2 \beta_{i-1} \beta_i + \beta_{i-1}^2)
+  (\beta_i^2 - 2 \beta_i \beta_{i-1} + \beta_{i-1}^2) 
\right) \\
&=& \frac{1}{12} b_i
\left( 12 \left(\frac{\beta_i  + \beta_{i-1}}{2}\right)^2
+  \left(\beta_i - \beta_{i-1}\right) ^2  \right) \\
&=& b_i m_i^2+\frac{1}{12} b_i^3.
\end{eqnarray*}

Summing over all intervals gives
\[ \frac13=\int_0^1x^2\,dx =\sum_{i=1}^s b_i m_i^2+\frac1{12}\sum_{i=1}^s b_i^3.
\]
Hence
\[ \frac13-\frac1{12}\sum_{i=1}^s b_i^3 = \sum_{i=1}^s b_i m_i^2
  \geq 
\sum_{i=1}^s b_i c_i^2 . 
\]
Since the nonnegative \(b_i\) sum to \(1\), at least one \(b_i>0\), and thus \(\sum_i b_i^3>0\). Therefore 
\[ \sum_{i=1}^s b_i c_i^2 \leq \frac13-\frac1{12}\sum_{i=1}^s b_i^3  < \frac13.\]
\end{proof}
\begin{remark} \label{RemarkAIuse}
This lemma was generated by  Claude, as a  revision of a ChatGPT proof
of a version of Theorem \ref{thm:barrier}.
%(distinct\_and\_nondistinct\_c\_proof\_revised.tex) }
We note that a previous version of our proof additionally assumed 
$0 \le c_1$ and $c_s \le 1$;  the Lean formalization (Section \ref{sec:lean}) revealed that these hypotheses  are never used.  This in turn simplified
the proof of Theorem \ref{thm:barrier} by eliminating the 
need to prove $c_i \in [0,1]$, which in turn allowed us to 
remove the   hypothesis $b^\top Ac = \tfrac16$.
This is a good illustration of the iterative process we used, and the continual
improvements that resulted from the intensive interaction between
Claude, ChatGPT, Gemini, Lean, and humans.
\end{remark}

 Finally, we are ready to prove Theorem \ref{thm:barrier}:
\begin{proof} [Proof of Theorem \ref{thm:barrier}:]
%\lajos{with this LaTeX command an extra . appears after "Proof"}
%\sigal{fixed!}
% We aim to prove that under the hypotheses of the theorem we have
% \begin{equation}\label{eq:quantbarrier}
%    b^\top (c \circ c) \;\le\; \frac13-\frac1{12}\sum_{i=1}^s b_i^3
%    \;\le\; \frac13 - \frac{1}{12\,s^2} \;<\; \frac13 .
% \end{equation}
The assumptions of the theorem are that $A, b \geq 0$ componentwise,
that the bushy tree conditions are satisfied: 
\[ b^\top \one =1, \; \;  b^\top c = \frac{1}{2},\; \; 
b^\top (c \circ c) = \frac{1}{3} ,\]
where $c = A \one$, and that the algebraic stability matrix $M$
is symmetric positive semi-definite.

Let $\chi_k\in\R^s$ be a vector with the first $k$ entries $1$ and the rest $0$.  The assumption on $M$ tells us that
\[ 0 \;\le\; \chi_k^\top  M \chi_k
   \;=\; 2\sum_{i,j\le k} b_i a_{ij}  - \Bigl(\sum_{i\le k}b_i\Bigr)^2
   \;=\; 2\sum_{i,j\le k} b_i a_{ij} - \beta_k^2 .\]
Since $A\ge0$, $b \ge 0$ and  $A \one = c$,
\[  \sum_{i=1}^k   \sum_{j = 1}^k b_i a_{ij}   \le  
\sum_{i=1}^k b_i \sum_{j=1}^s a_{ij}  = \sum_{i=1}^k b_i c_i   .
\]
Plugging this back into the inequality above we obtain 
\[ 0 \leq  2\sum_{i,j\le k} b_i a_{ij} - \beta_k^2   \leq 2 \sum_{i=1}^k b_i c_i  - \beta_k^2 .\]
This gives $\sum_{i=1}^k b_i c_i  \ge \tfrac12 \beta_k^2$ for every $k$.
Now the hypotheses of Lemma
\ref{lem:partialsums} are satisfied and so we establish that
$ b^\top (c \circ c)  < \tfrac13,$
contradicting the third order bushy tree condition \eqref{eq:RKOC3a}.
% \[ b^\top (c \circ c) \;\le\; \frac13-\frac1{12}\sum_{i=1}^s b_i^3.\]
% We know that $b_i \geq 0$ and that $\sum_i b_i =1 $, so there must be at least on
% value of $b_i$ that is strictly positive. This means that
% $ \sum_{i=1}^s b_i^3 > 0,$
% so that 

\end{proof}

\begin{remark}
    We can better understand by {\em how much} the third order condition is not 
    met. Observe that  the power-mean inequality gives
    $ \sum_i b_i^3 \ge \bigl(\sum_i b_i\bigr)^3/s^2 = 1/s^2,$
    so we have
    \begin{equation}\label{eq:quantbarrier}
   b^\top (c \circ c) \;\le\; \frac13-\frac1{12}\sum_{i=1}^s b_i^3
   \;\le\; \frac13 - \frac{1}{12\,s^2} \;<\; \frac13 .
\end{equation}
\end{remark}

\section{Machine verification}\label{sec:lean}

\begin{table}[ht]
\centering
\begin{threeparttable}
\caption{Correspondence between results in this paper and their Lean formalizations
(identifiers are in the namespace \texttt{SSP.MainRevised}).}
\label{tab:lean-correspondence}
\begin{tabular}{@{}lll@{}}
\toprule
Result & Lean identifier & Status \\
\midrule
\resref{Theorem}{thm:SSPbound}
  & \leanproof{ssp\_bound\_of\_conditions}{MainRevised}{553}
  & Verified \\
\resref{Lemma}{lem:wcs}
  & \leanproof{weighted\_cauchy\_schwarz}{MainRevised}{83}
  & Verified \\
\resref{Lemma}{lem:substochastic}
  & \leanproof{substochastic\_rowSum\_le}{MainRevised}{91} & Verified\\
  &  \leanproof{substochastic\_entry\_mem}{MainRevised}{105}& Verified\\
  &  \leanproof{substochastic\_neumann}{MainRevised}{622}  & Verified\\%\tnote{4} \\
\resref{Lemma}{lem:neumann2}
  & \leanproof{neumann\_square\_series}{MainRevised}{650}
  & Verified \\
\resref{Lemma}{lem:AS}
  & \leanproof{substochastic\_sequence\_inequality}{MainRevised}{660}
  & Verified \\ %\tnote{6} \\
\resref{Lemma}{lem:componentbound}
  & \leanproof{componentwise\_bound}{MainRevised}{229}
  & Verified \\
\resref{Theorem}{thm:alt}
  & \leanproof{ssp\_bound\_of\_invMat}{MainRevised}{247}
  & Verified \\
\resref{Theorem}{thm:AldousShepp}
  & \leanproof{aldous\_shepp}{MainRevised}{671}
  & Verified \\%\tnote{11} \\
\resref{Theorem}{thm:OsogamiHarcholBalter}
  & \leanproof{osogami\_harchol\_balter}{MainRevised}{701}
  & Verified \\
\resref{Theorem}{thm:acyclic-hierarchy}
  & \leanproof{acyclic\_hierarchy\_ratio}{MainRevised}{728},
    \leanproof{acyclic\_hierarchy}{MainRevised}{713}
  & Verified \\ %\tnote{13} \\
\resref{Theorem}{thm:cyclic-hierarchy}
  & \leanproof{cyclic\_hierarchy}{MainRevised}{685}
  & Verified \\ %\tnote{14} \\
\resref{Lemma}{lem:alg}
  & \leanproof{algebraic\_inequality}{MainRevised}{744}
  & Verified \\ %\tnote{15} \\
\resref{Theorem}{thm:SSPboundP3}
  & \leanproof{order3\_bound\_of\_invMat}{MainRevised}{753}
  & Verified \\ %\tnote{16} \\
\resref{Theorem}{thm:barrier}
  & \leanproof{barrier\_unsorted}{MainRevised}{823}
  & Verified \\ %\tnote{17} \\
\resref{Lemma}{lem:partialsums}
  & \leanproof{partial\_sums\_strict}{MainRevised}{804}
  & Verified \\
\bottomrule
\end{tabular}
% \begin{tablenotes}\footnotesize
% \item[4] Proved in these parts, including spectral radius, convergent Neumann series and inverse nonnegativity.
% \item[6] Proved in series form.
% \item[11] Proved for concrete transient initial laws with at most $s$ states, including an atom at zero.
% \item[13] Proved in triangular uniformized-density form, including positive-moment ratios, at most $s$ states, and an atom at zero; the division-free form also covers zero moments.
% \item[14] Proved in ratio form for concrete transient initial laws with at most $s$ states, including an atom at zero.
% \item[15] Proved, including $w=0$ and $\beta=0$; Lean division by zero gives the stated convention.
% \item[16] Proved directly from inverse SSP conditions and transformed order-three identities.
% \item[17] Proved, including a permutation to sort stages; no sorted-stage hypothesis.
% \end{tablenotes}
\end{threeparttable}
\end{table}

The elementary content of this paper has been formalized and verified in Lean~4 with
Mathlib (repository directory \texttt{SSPLean/}; all results compile with no
\texttt{sorry} and depend only on Lean's three standard axioms).  
\resref{Table}{tab:lean-correspondence} lists
the correspondences.

% \begin{center}
% \small
% \begin{tabular}{lll}
% \toprule
% result & Lean identifier & status\\
% \midrule
% Lemma \ref{lem:partialsums} (generalized) & \texttt{partial\_sum\_lemma} & verified\\
% Theorem \ref{thm:barrier} & \texttt{barrier\_bound}, \texttt{no\_order3\_algebraically\_stable} & verified\\
% quantitative barrier \eqref{eq:quantbarrier} & \texttt{quantitative\_barrier} & verified\\
% Lemma \ref{lem:wcs} & \texttt{weighted\_cauchy\_schwarz} & verified\\
% Lemma \ref{lem:substochastic} (finite part) & \texttt{entry\_le\_one\_of\_substochastic} & verified\\
% Theorem \ref{thm:SSPbound}  & \texttt{ssp\_bound\_of\_conditions} & verified$^{*}$\\
% Lemma \ref{lem:componentbound} & \texttt{ComponentBound} & verified$^{*}$\\
% Theorem \ref{thm:order3} (reduction) & 
% \texttt{order3\_bound\_unconditional}
% & verified\\
% Conjecture \ref{conj:hierarchy}, $k=2$ & \texttt{hierarchy\_k2} & verified\\
% \bottomrule
% \end{tabular}

% \smallskip
% {\footnotesize $^{*}$\,verified using a proof which avoids using the
% Aldous--Shepp theorem in (Lemma \ref{lem:componentbound})}
% \end{center}
Lemma \ref{lem:componentbound} was formalized using an explicit weighted path interpretation which avoided any use of 
the Aldous--Shepp theorem or the Perron-Frobenius theorem 
so the  $r \le 2s$ result is machine-checked end
to end.  The formalization process found the simplifications described in  Remark \ref{RemarkAIuse}.

\section*{On our collaboration}

Our collaboration started at the SIAM CSE meeting in 2025 (Dallas, TX) 
after Sigal Gottlieb (SG) gave a talk on strong stability preserving 
time-discretizations.  Lior Horesh (LH) noted that AI tools could be 
valuable for proving SSP bounds and barriers, and obtaining optimal methods. 
LH and SG arranged a meeting with colleagues David Ketcheson (DK), 
Zachary J. Grant (ZJG), Avi Shinnar (AS), and Barry Trager (BT). 
At this meeting, SG, ZJG, and DK proposed the study of some facts 
that were believed to be true about Runge--Kutta methods, 
but were not proved (including the 2s conjecture, Conjecture \ref{conj:SSPbound}). 
%\lajos{was one of the facts the impossibility of the 3rd-order non-negative A-stability conjecture which we removed from this manuscript?} \david{Yes.} 
This was the beginning of the collaboration. 
Later, LL joined our team and used his expertise to validate
some of our results, prove other results, and explain  their significance.

In the initial phase we needed to understand how to formulate the 
problems in a way amenable to AI tools. Once we learned to do this  
we used different AI tools to prove certain results, 
and cross-checked it with other AI tools to identify gaps, 
inconsistencies, and misstatements. In many cases, the AI tools 
gave us wrong proofs or mistaken implications; 
these were often persistent and the AI required careful 
checking and firm instructions to avoid these. 
The interaction between AI tools and the mathematical
researchers was iterative.  Once we were satisfied that a 
result was correct, we processed it  into our own mathematical 
language and formalized it in Lean to verify the results.
This often led to further refinements and more clarity as to the structure
of the proof and which hypotheses were actually needed.
For example, the Lean formalization (Section \ref{sec:lean}) found that one 
hypothesis of the key barrier lemma was unnecessary, 
which in turn removed an entire step of the barrier proof and
one order condition from its hypotheses.

\appendix 
\section{An elementary algebraic proof of Lemma \ref{lem:componentbound}}\label{sec:AppendixA}
In this section, we present a self-contained elementary linear-algebraic proof of the entrywise substochastic inequality \eqref{eq:componentbound}. This proof operates purely via matrix identities, sum-of-squares decompositions, and finite-dimensional Cauchy--Schwarz inequalities, avoiding all references to continuous-time Markov chains, martingales, phase-type distributions, or graph theory. The proof was derived through a human--AI collaboration between LL and Gemini and then verified 
independently by ChatGPT  and re-written by SG.  

First, we point the reader to Remark 
\ref{rem:cjnondecreasing}, that explains that
for a Runge--Kutta method given by $A,b$, we can always
permute the rows and columns of  $A$ and the entries of  $b$
so that the abscissas $c_i$ are non-decreasing.
We use this in the following where we  assume 
\[ c_1 \le c_2 \le \cdots \le c_s \]
without loss of generality. For this reason, we restate the lemma 
with the assumption of non-decreasing $c_i$, without loss of generality.

\smallskip

\noindent {\bf Restatement of Lemma \ref{lem:componentbound}:} Let $s$ be a positive integer, and let $Q \in \mathbb{R}^{s \times s}$ be an entrywise non-negative matrix ($Q \ge 0$) that is substochastic ($Q\mathbbm{1} \le \mathbbm{1}$) with spectral radius $\rho(Q) < 1$. 
Let $N := (I - Q)^{-1} = \sum_{k=0}^{\infty} Q^k$, and define the fundamental vector 
$x := N\one$, where $0\le x_1\le x_2\le\cdots\le x_s$.
Then, for any entry index $i \in \{1, 2, \dots, s\}$, the following inequality holds
\begin{equation}\label{targetstatement}
x_i^2 \le \frac{2s}{s+1} (Nx)_i.
\end{equation}

\smallskip

\begin{proof}[An algebraic proof of Lemma \ref{lem:componentbound}]
Set \(x=N \one\), 
so that
\[  N = (I- Q)^{-1} \; \;  \Rightarrow \; \; Q x = x-\one. \]
Also, 
\[ N=\sum_{\ell=0}^{\infty}Q^\ell\ge0
\; \;  \Rightarrow \; \; x \ge 0.
\]
Moreover, recall that $N = I + rA$ so that the assumption
that $c_1\le c_2\le\cdots\le c_s$ coupled with 
$ x = N \one = \one + r c  $ means that
\[  
0\le x_1\le x_2\le\cdots\le x_s .
\]

Set \(x_0=0\), and define
\[ \Delta_m=x_m-x_{m-1}\ge0, \qquad \mbox{and} \qquad
h_k=\sum_{m=1}^k \Delta_m^2 .
\]
Since the differences are nonnegative we can conclude that
\[
h_k=\sum_{m=1}^k\Delta_m^2
\le\left(\sum_{m=1}^k\Delta_m\right)^2=x_k^2.
\]
Moreover, when \(k>j\),
\[
h_k-h_j
=\sum_{m=j+1}^k\Delta_m^2
\le\left(\sum_{m=j+1}^k\Delta_m\right)^2
=(x_k-x_j)^2.
\]
When \(k\le j\), the left side is nonpositive, and this inequality is still true, so that 
\[ h_k-h_j \le (x_k-x_j)^2, \qquad \forall j, k. \]

Define the vector $h$ comprised of the elements $h_k$, and let 
\(d= N^{-1} \one = \one -Q \one\ge0\) so that
\( d_k = 1 - \sum_{j=1}^s Q_{kj} .\)
The  inequalities above give 
\[ \begin{aligned}
((I-Q)h)_k & =h_k-\sum_{j=1}^s Q_{kj} h_j
=\left(1-\sum_{j=1}^sQ_{kj}\right)h_k
  +\sum_{j=1}^sQ_{kj}(h_k-h_j)
=d_kh_k+\sum_{j=1}^sQ_{kj}(h_k-h_j) \\
&\le\sum_{j=1}^sQ_{kj}(x_k-x_j)^2+d_kx_k^2 
=\left(\sum_{j=1}^sQ_{kj}+d_k\right)x_k^2
-2x_k\sum_{j=1}^sQ_{kj}x_j +\sum_{j=1}^sQ_{kj}x_j^2 \\
&=x_k^2-2x_k(Qx)_k+\sum_{j=1}^sQ_{kj}x_j^2
= x_k^2-2x_k(x_k -1)+\sum_{j=1}^sQ_{kj}x_j^2 
=2x_k-x_k^2+\sum_{j=1}^sQ_{kj}x_j^2.
\end{aligned}
\]
Since this is true for each element $k$, we have 
\[ (I-Q) h \le 2 x - (I-Q)(x\circ x).
\]
Multiplying by \(N\ge0\) gives 
\[
h \le 2 N x - x\circ x
\, \, \Rightarrow \; \; 
h_i \le 2 (N x)_i - x_i^2.
\]
Finally, we use this and Cauchy–Schwarz  to obtain 
\[ x_i^2
=\left(\sum_{m=1}^i\Delta_m\right)^2
\le i\sum_{m=1}^i\Delta_m^2
\le s \, h_i
\le s \, \bigl(2(Nx)_i-x_i^2\bigr).
\]
Rearranging and using \(x=N \one\) we get the result:
\[
(s+1)\bigl((N \one)_i\bigr)^2
\le 2s(N^2 \one)_i .
\]
\end{proof}

\section{Proof of a simple Lemma} \label{app:LemAlg}
In this section we present the proof of Lemma \ref{lem:alg}, we restate it here:

\noindent{\bf Lemma \ref{lem:alg}}
{\em Let $n\ge 2$, $w\in[0,1)$, $u\ge0$, and
$\beta\ge\frac{n}{n-1}\,w^2$; we define $\beta = 0$  when $w=0$,
and in particular $\beta^2/w = 0$  when $w=0$. Put
\[ m = 2(1-w)+u+\beta . \]
Then
\begin{equation}\label{eq:alg}
  3(1-w)\,m + 3u + \frac{n+1}{n} \left(\frac{u^2}{1-w}+\frac{\beta^2}{w}\right)
  \;\ge\; \frac{n+2}{n+1}\,m^2 .
\end{equation}
}

\begin{proof}
\emph{(a) Eliminate the fractions by Cauchy--Schwarz.} The Engel form of
Cauchy--Schwarz states
$\frac{X^2}{\alpha}+\frac{Y^2}{\gamma}\ge\frac{(X+Y)^2}{\alpha+\gamma}$
for $\alpha,\gamma>0$. Apply it with $X=u+(1-w)$, $Y=\beta$,
$\alpha=1-w$, $\gamma=w$. Since $\alpha+\gamma=1$ and
$X+Y=m-(1-w)$,
\[
  \frac{\bigl(u+(1-w)\bigr)^2}{1-w}+\frac{\beta^2}{w}
  \;\ge\; \bigl(m-(1-w)\bigr)^2 .
\]
When $w=0$ this still holds, since then $\beta=0$ and both sides equal
$(u+1)^2$. The first term on the left equals
$\frac{u^2}{1-w}+2u+(1-w)$, so
\begin{equation}\label{eq:CS}
  \frac{u^2}{1-w}+\frac{\beta^2}{w}
  \;\ge\; \bigl(m-(1-w)\bigr)^2-2u-(1-w).
\end{equation}

\emph{(b) Reduce to a polynomial.} Insert \eqref{eq:CS} into the
left-hand side of \eqref{eq:alg}, subtract $\frac{n+2}{n+1}m^2$, and
expand. Using $\frac{n+1}{n} -\frac{n+2}{n+1}=\frac{1}{n(n+1)}$ and
$3-2\frac{n+1}{n}=\frac{n-2}{n}$, this gives
\[
  \text{LHS of \eqref{eq:alg}}-\frac{n+2}{n+1}m^2
  \;\ge\; \Phi(m,u,w),
\]
where
\[
  \Phi(m,u,w)=\frac{m^2}{n(n+1)}
  +\frac{n-2}{n}(1-w)\,m
  +\frac{n+1}{n}(1-w)^2
  +\frac{n-2}{n}\,u
  -\frac{n+1}{n}(1-w).
\]
So it suffices to show $\Phi(m,u,w)\ge0$.

\emph{(c) Reduce to the extreme case $u=0$, $\beta=\frac{n}{n-1}w^2$.}
For $n\ge2$ and $m>0$,
\[
  \frac{\partial\Phi}{\partial m}
  =\frac{2m}{n(n+1)}+\frac{n-2}{n}(1-w)\;\ge\;0,
  \qquad
  \frac{\partial\Phi}{\partial u}=\frac{n-2}{n}\;\ge\;0 .
\]
So $\Phi$ is nondecreasing in $m$ on $(0,\infty)$ and in $u$. By the
hypotheses on $u$ and $\beta$,
\[
  m \;\ge\; M_0(w):=2(1-w)+\frac{n}{n-1}\,w^2
  =\frac{n+1}{n}+\frac{n}{n-1}\Bigl(w-\frac{n-1}{n}\Bigr)^2\;>\;0 .
\]
Hence $\Phi(m,u,w)\ge\Phi\bigl(M_0(w),0,w\bigr)=:\ell(w)$.

\emph{(d) Factor $\ell$.} Direct expansion shows that
\begin{equation}\label{eq:factor}
  \ell(w)
  =\frac{1}{n(n+1)}\,
   \bigl(nw-(n-1)\bigr)^2
   \left(\frac{w^2}{(n-1)^2}-w+\frac{2n}{n-1}\right).
\end{equation}
The first factor is a square. For $0\le w\le 1$ the second factor
satisfies
\[
  \frac{w^2}{(n-1)^2}-w+\frac{2n}{n-1}
  \;\ge\; -1+\frac{2n}{n-1}
  \;=\;\frac{n+1}{n-1}\;>\;0 .
\]
Therefore $\ell(w)\ge0$, so $\Phi(m,u,w)\ge0$, which proves
\eqref{eq:alg}.
\end{proof}

\bibliographystyle{plain}
\bibliography{refs}

\begin{thebibliography}{10}

\bibitem{AldousShepp1987}
D.~Aldous and L.~Shepp.
\newblock The least variable phase type distribution is {E}rlang.
\newblock {\em Stochastic Models}, 3(3):467--473, 1987.

\bibitem{Asmussen2003}
S.~Asmussen.
\newblock {\em Applied Probability and Queues}.
\newblock Springer, second edition, 2003.

\bibitem{bolley1978}
C.~Bolley and M.~Crouzeix.
\newblock {Conservation de la positivit\'{e} lors de la discr\'{e}tisation des
  probl\'{e}mes d'\'{e}volution paraboliques}.
\newblock {\em R.A.I.R.O. Analyse Num{\'{e}}rique}, 12(3):237--245, 1978.

\bibitem{lean4}
L.~de~Moura and S.~Ullrich.
\newblock The {L}ean 4 theorem prover and programming language.
\newblock In {\em CADE-28: The 28th International Conference on Automated
  Deduction}. Springer-Verlag, 2021.

\bibitem{FerracinaSpijker2004}
L.~Ferracina and M.~N. Spijker.
\newblock Stepsize restrictions for the total-variation-diminishing property in
  general {R}unge--{K}utta methods.
\newblock {\em SIAM Journal on Numerical Analysis}, 42:1073--1093, 2004.

\bibitem{FerracinaSpijker2008}
L.~Ferracina and M.~N. Spijker.
\newblock Strong stability of singly-diagonally-implicit {R}unge--{K}utta
  methods.
\newblock {\em Applied Numerical Mathematics}, 58(11):1675--1686, 2008.

\bibitem{GKSbook}
S.~Gottlieb, D.~I. Ketcheson, and C.-W. Shu.
\newblock {\em Strong Stability Preserving Runge--Kutta and Multistep Time
  Discretizations}.
\newblock World Scientific, 2011.

\bibitem{HairerWanner}
E.~Hairer and G.~Wanner.
\newblock {\em Solving Ordinary Differential Equations II: Stiff and
  Differential-Algebraic Problems}.
\newblock Springer, second edition, 1996.

\bibitem{KemenySnell1976}
J.~G. Kemeny and J.~L. Snell.
\newblock {\em Finite Markov Chains}.
\newblock Springer, 1976.

\bibitem{RKopt}
D.~Ketcheson, M.~Parsani, Z.~Grant, A.~Ahmadia, and H.~Ranocha.
\newblock {RK-Opt}: A package for the design of numerical {ODE} solvers.
\newblock {\em Journal of Open Source Software}, 5(54):2514, 2020.

\bibitem{KMG09}
D.~I. Ketcheson, C.~B. Macdonald, and S.~Gottlieb.
\newblock Optimal implicit strong stability preserving {R}unge--{K}utta
  methods.
\newblock {\em Applied Numerical Mathematics}, 59:373--392, 2009.

\bibitem{KocsisNemeth}
T.~A. Kocsis and A.~N{\'e}meth.
\newblock Optimal second order diagonally implicit {SSP} {R}unge--{K}utta
  methods, 2014.

\bibitem{Kraaijevanger1991}
J.~F. B.~M. Kraaijevanger.
\newblock Contractivity of {R}unge--{K}utta methods.
\newblock {\em BIT}, 31:482--528, 1991.

\bibitem{LatoucheRamaswami1999}
G.~Latouche and V.~Ramaswami.
\newblock {\em Introduction to Matrix Analytic Methods in Stochastic
  Modelling}.
\newblock SIAM, 1999.

\bibitem{2014_rational_absmon}
L.~L{\'o}czi and D.~I. Ketcheson.
\newblock Rational functions with maximal radius of absolute monotonicity.
\newblock {\em LMS Journal of Computation and Mathematics}, 17(1):159--205,
  2014.

\bibitem{Neuts1981}
M.~F. Neuts.
\newblock {\em Matrix-Geometric Solutions in Stochastic Models: An Algorithmic
  Approach}.
\newblock Johns Hopkins University Press, 1981.

\bibitem{Osogami2005}
T.~Osogami.
\newblock {\em Analysis of Multi-server Systems via Dimensionality Reduction of
  Markov Chains}.
\newblock PhD thesis, Carnegie Mellon University, June 2005.

\bibitem{OsogamiHarcholBalter2006}
T.~Osogami and M.~Harchol-Balter.
\newblock Closed form solutions for mapping general distributions to
  quasi-minimal {PH} distributions.
\newblock {\em Performance Evaluation}, 63(6):524--552, 2006.

\bibitem{ShuOsher1988}
C.-W. Shu and S.~Osher.
\newblock Efficient implementation of essentially non-oscillatory
  shock-capturing schemes.
\newblock {\em Journal of Computational Physics}, 77(2):439--471, 1988.

\bibitem{ShuOsher1989}
C.-W. Shu and S.~Osher.
\newblock Efficient implementation of essentially non-oscillatory
  shock-capturing schemes, {II}.
\newblock {\em Journal of Computational Physics}, 83(1):32--78, 1989.

\bibitem{spijker1983}
M.~N. Spijker.
\newblock {Contractivity in the numerical solution of initial value problems}.
\newblock {\em Numerische Mathematik}, 42:271--290, 1983.

\bibitem{Tao2026AI}
T.~Tao.
\newblock Mathematics in the age of {AI}, 2026.

\bibitem{Thurston}
W.~P. Thurston.
\newblock On proof and progress in mathematics.
\newblock {\em Bull. Amer. Math. Soc.}, 30(2):161--177, 1994.

\bibitem{vandegriend1986}
J.~A. van~de Griend and J.~F. B.~M. Kraaijevanger.
\newblock Absolute monotonicity of rational functions occurring in the
  numerical solution of initial value problems.
\newblock {\em Numerische Mathematik}, 49:413--424, 1986.

\end{thebibliography}

\end{document}